\documentclass{article}
\usepackage[utf8]{inputenc}
\usepackage{amsthm}
\usepackage{amsfonts}         
\usepackage{amsmath}
\usepackage{amssymb}
\usepackage{graphicx}
\usepackage{color}
\usepackage{float}
\usepackage{hyperref}
\usepackage{doi}   
\usepackage{url}   

\usepackage{longtable}
\usepackage[pagewise]{lineno}

\newcommand{\R}{\mathbb{R}}

\newcommand{\argmin}{\mathop{\mathrm{arg\,min}}}

\title{Dual-grid parameter choice method with application to image deblurring}

\author{\underline{Markus Juvonen}$^1$\footnote{Corresponding author: Markus Juvonen, email: markus.juvonen@helsinki.fi}, Bj\o{}rn Jensen$^1$, Ilmari Pohjola$^1$, \\
Yiqiu Dong$^2$ and Samuli Siltanen$^1$}

\date{%
    \normalsize\textit{$^1$University of Helsinki, Finland}\\%
    \normalsize\textit{$^2$Technical University of Denmark, Denmark}\\[2ex]%
    April 2025}

\begin{document}

\maketitle

\begin{abstract}
Variational regularization of ill-posed inverse problems is based on minimizing the sum of a data fidelity term and a regularization term. The balance between them is tuned using a positive regularization parameter, whose automatic choice remains an open question in general. A novel approach for parameter choice is introduced, based on the use of two slightly different computational models for the same inverse problem. Small parameter values should give two very different reconstructions due to amplification of noise. Large parameter values lead to two identical but trivial reconstructions. The optimal parameter is chosen between the extremes by matching image similarity of the two reconstructions with a pre-defined value. The efficacy of the new method is demonstrated by image deblurring using measured data and two different regularizers. 
\end{abstract}


\section{Introduction}

Variational regularization is a fundamental solution technique for ill-posed inverse problems. The basic idea is to minimize a penalty functional with two terms: a data fidelity term and a regularization term, the former ensuring solutions remain close to the measured data and the latter modeling {\it a priori} information about the unknown \cite{engl1996regularization}. The balance between the two terms is tuned using a positive regularization parameter, whose automatic choice remains an open question, in general. We introduce a novel approach for parameter choice based on the use of two slightly different computational models. 

Consider a forward problem model
\begin{equation}\label{basicmodel}
m = \mathcal{A}(f) + \varepsilon,
\end{equation}
where $m$ is the observed measurement data, $f$ is the underlying unknown true target, $\varepsilon$ is measurement noise and $\mathcal{A}$ is the measurement model. Our formulation of variational regularization is 
\begin{equation}\label{regul_linear2}
\argmin_{f} \|{m} - Af\|_2^2 + \alpha \mathcal{R}({f})
\end{equation}
where $\alpha>0$ is the regularization parameter and $A$ is a matrix. The data fidelity term $\|{m} - Af\|_2^2$ penalizes vectors $f$ for not matching the measured data, and the regularizer $\mathcal{R}$ should be picked or designed so that the value $\mathcal{R}(f)>0$ is larger for any $f$ that are not expected in light of any {\it a priori} information. We focus on choosing $\alpha$ automatically in an optimal way.

Many methods have been suggested in the literature: discrepancy principles \cite{anzengruber2009morozov,toma2015iterative,wen2011parameter}, cross-validation \cite{lukas2008strong, wen2018using}, the L-curve method \cite{hansen1992analysis}, sparsity matching \cite{hamalainen2013sparse, meaney2024image,purisha2018automatic,purisha2017controlled}, multigrid method \cite{niinimaki2016multiresolution}, residual whiteness principle \cite{pragliola2023admm}, as well as methods tailored to select parameters for total variation regularization \cite{chen2014automatic, clason2010duality, davoli2024dyadic, dong2011automated, kindermann2014numerical}. However, while some of the above methods might work for a given application, they may not be useful for slightly different applications.
A complementary line of work is bilevel parameter learning, where the regularization weight(s) are learned from examples by minimizing an upper-level task loss while the lower level solves the variational model \cite{arridge2019solving,delosreyes2017bilevel,kunisch2013bilevel}. 
In practice, automatic parameter choice remains a challenge.

We introduce a new parameter choice principle, similar in spirit to the multi-resolution idea in \cite{niinimaki2016multiresolution}, but simpler and more general. Our approach is a finite-dimensional heuristic intended for realistic use rather than worst-case guarantees. While the Bakushinskii veto \cite{bakushinskii1984remarks} rules out deterministic, infinite-dimensional noise-level-free rules in the worst case, practical success of heuristics is well documented, and the veto can be overcome or avoided in such settings \cite{bauer2008regularization,harrach2020beyond,neubauer2008convergence}.

Our idea is to use two forward models, both written for the same data $m$ but using different computational grids to represent the unknown $f$. We then solve the same problem for both grids, for multiple parameter values $\alpha$, and compare the structural similarity of the two solutions using SSIM \cite{SSIM}. If the parameter is too small, we expect to have a low similarity between the solutions, as noise is amplified in an ill-posed reconstruction process. However, the similarity should grow with increasing $\alpha$ when the problem stabilizes. We can set a suitable target threshold for the SSIM value and find a good parameter value for the reconstructions as we reach this threshold.

In this paper, we study the new method in the context of image deblurring. 
Our approach is independent of the regularizer, so we test it with two different ones: Tikhonov regularization \cite{hansen2006deblurring,Tik63} and total variation (TV) regularization \cite{doi:10.1137/040605412, RudinOsher94,VogelTV98, doi:10.1137/080724265}.
We test our method with both simulated and real data and compare the new algorithm with the classical discrepancy principle and a bilevel optimization method.




\section{Materials and methods} 

\subsection{Dual-grid model for image deblurring} 

Let us explain our approach in more detail using two-dimensional convolution 
\begin{equation}\label{measmodel1}
  m(x) = (p\ast f)(x) = \int_{\R^2} f(x-y)p(y) dy 
\end{equation}
as the model problem. In (\ref{measmodel1}), measured data is again denoted by $m$, the unknown is $f$, and $p$ is a point spread function. Finding $f$ when $p$ is given and $m$ measured is the inverse problem called deconvolution. 

Let $s\in\R^2$ be a displacement vector called the {\it shift}. Define the shifted image $f_s$ and shifted point spread function $p_{-s}$ as follows
$$
  f_s(x) = f(x+s), \qquad    p_{-s}(x) = p(x-s).
$$
Substituting $\lambda=y+s$ to (\ref{measmodel1}) we get the shift identity
\begin{eqnarray}\label{shifteq}
  (p\ast f)(x) 
  = \int_{\R^2} f_s(x-\lambda)p_{-s}(\lambda) d\lambda  =   (p_{-s}\ast f_s)(x).
\end{eqnarray}
Now consider another forward model using the shifted point spread function:
\begin{equation}\label{measmodel2}
  m = p_{-s}\ast g.
\end{equation}
 Assume (unrealistically) that we can solve the inverse problems (\ref{measmodel1}) and  (\ref{measmodel2}) perfectly and that the measurement in both cases is the same $m$. Then (\ref{shifteq}) implies  that the recovered images are related by a shift: $g = f_s$.

Let's move to a more realistic  model. 
Given a blurred pixel image $\mathbf{m}\in \R^{n{\times}n}$, we use a discrete computational model with periodic boundary conditions:
\begin{equation}\label{realisticmodel}
  \mathbf{m} 
  = 
  \mathbf{p}\star \mathbf{f} +\mathbf{e} 
  =
    \mathcal{F}^{-1}\big(\mathcal{F}(\mathbf{p})\cdot \mathcal{F}(\mathbf{f})\big) +\mathbf{e},
\end{equation}
where $\mathcal{F}$ is the discrete Fourier transform, $\mathbf{f}\in \R^{n{\times}n}$ is the sharp image defined on the same pixel grid where $ \mathbf{m}$ is given (camera sensor pixels),  $\mathbf{p}\in \R^{n{\times}n}$ is the discrete PSF, and $\mathbf{e}\in \R^{n{\times}n}$ models random measurement noise. Also, ``$\cdot$'' in (\ref{realisticmodel}) means the pixel-wise Hadamard product.

Denote the length of the side of a square pixel in our grid by $h>0$. For the rest of this study, we fix the shift vector to be
$$
  s = [h/2,h/2]^T \in\R^2.
$$
Define ``half-pixel shift kernels'' related to the shift vector $s$:
$$
\mathbf{d}_s = 
\left[\!\!
\begin{array}{ccccc}
0 &0 &0 &0 &0 \\
0 &0 &0 &0 &0 \\
0 &1/4 &\mathbf{1/4} &0 &0 \\
0 &1/4 &1/4 &0 &0 \\
0 &0 &0 &0 &0 
\end{array}
\!\!\right],
\qquad
\mathbf{d}_{-s} = 
\left[\!\!
\begin{array}{ccccc}
0 &0 &0 &0 &0 \\
0 &0 &1/4 &1/4 &0 \\
0 &0 &\mathbf{1/4} &1/4 &0 \\
0 &0 &0 &0 &0 \\
0 &0 &0 &0 &0 
\end{array}
\!\!\right],
$$
where there need to be an appropriate amount of zero elements according to image size. The central elements of the kernels are indicated as bold. 

Consider two computational models: the camera grid model
\begin{eqnarray} \label{discmeasmodel1}
  \mathbf{m} =   \mathbf{p}\star \mathbf{f} +\mathbf{e},
\end{eqnarray}
and the shifted model
\begin{eqnarray} \label{discmeasmodel2}
  \mathbf{m} = \mathbf{p}_{-s} \star \mathbf{g} +\mathbf{e},
\end{eqnarray}
where  $\mathbf{g}\in \R^{n{\times}n}$ is the desired sharp image {\it on a shifted pixel grid} and 
$\mathbf{p}_{-s} = \mathbf{p}\star \mathbf{d}_{-s} $ is the shifted PSF. 
Note that the data $ \mathbf{m}$ is the same in both models (\ref{discmeasmodel1}) and  (\ref{discmeasmodel2}).

Solve the inverse problems (\ref{discmeasmodel1}) and  (\ref{discmeasmodel2}) by regularization with parameter $\alpha>0$ and regularizer $\mathcal{R}:\R^{n{\times}n}\rightarrow\R^+$, yielding reconstructions 
\begin{eqnarray}\label{genrealregularrecon}
{\mathbf{f}}^{(\alpha)} &=& \argmin_{\mathbf{f}\in \R^{n{\times}n}}\|\mathbf{m} - \mathbf{p}\star \mathbf{f}\|^2_F + \alpha \mathcal{R}(\mathbf{f}),\\
\label{genrealshiftedrecon}
{\mathbf{g}}^{(\alpha)} &=& \argmin_{\mathbf{g}\in \R^{n{\times}n}}\|\mathbf{m} - \mathbf{p}_{-s} \star \mathbf{g}\|^2_F + \alpha \mathcal{R}(\mathbf{g}),
\end{eqnarray}
\textcolor{black}{where $\|\cdot\|_F$ denotes the Frobenius norm.}
The above discussion of the continuous problem suggests that with any large enough $\alpha>0$ we should have 
$
  {\mathbf{g}}^{(\alpha)} \approx {\mathbf{f}}^{(\alpha)}\star \mathbf{d}_s.
$
Consequently the structural similarity of the two images should be high: SSIM$({\mathbf{g}}^{(\alpha)} ,{\mathbf{f}}^{(\alpha)}\star \mathbf{d}_s)\approx 1$.

Furthermore, convolving a properly regularized image ${\mathbf{f}}^{(\alpha)}$ by $\mathbf{d}_s$ does not change the image much under SSIM.  Therefore we get 
$$
  \mbox{SSIM}({\mathbf{g}}^{(\alpha)} ,{\mathbf{f}}^{(\alpha)})\approx 1.
$$

However, with too small parameter $\alpha>0$ the noise in the measurements is uncontrollably amplified in the reconstruction process, resulting in erratic pixel values in ${\mathbf{g}}^{(\alpha)}$ and ${\mathbf{f}}^{(\alpha)}$. The noise components in the two reconstructions are expected to be different because of \textcolor{black}{the randomness behaviour.}

We arrive at our research hypothesis:
\begin{itemize}
\item If the regularization parameter $\alpha>0$ is too small, then $\mbox{SSIM}({\mathbf{g}}^{(\alpha)},{\mathbf{f}}^{(\alpha)})$ is \textcolor{black}{much smaller than 1,} 
indicating that ${\mathbf{f}}^{(\alpha)}$ and ${\mathbf{g}}^{(\alpha)}$ are two very different images due to amplification of noise in the unstable inversion process.
\item If $\alpha>0$ is large enough, then $\mbox{SSIM}({\mathbf{g}}^{(\alpha)},{\mathbf{f}}^{(\alpha)})$ is close to 1. Then ${\mathbf{g}}^{(\alpha)}$ and ${\mathbf{f}}^{(\alpha)}$ are approximately equal, and so regularization is not too weak.
\end{itemize}
Our proposed strategy is to pick the smallest $\alpha$ for which $\mbox{SSIM}({\mathbf{g}}^{(\alpha)},{\mathbf{f}}^{(\alpha)})\geq T$, where \textcolor{black}{$T\in(0,1)$ and close to 1} is a user-defined lower bound for image quality.

To back up the hypothesis, we performed a numerical experiment reported below in Section \ref{sec:hypo_backup}.

\subsection{Tikhonov regularization}

The Tikhonov regularized solution of \eqref{genrealregularrecon} and \eqref{genrealshiftedrecon} is defined on the two grids by: 

\begin{eqnarray}\label{def:discreteTikhregu}
{\mathbf{f}}^{(\alpha)} &=& \argmin_{\mathbf{f}\in \R^{n{\times}n}}\|\mathbf{m} - \mathbf{p}\star \mathbf{f}\|^2_F + \alpha \|\mathbf{f}\|^2_F, \\
\label{def:discreteTikhshift}
{\mathbf{g}}^{(\alpha)} &=& \argmin_{\mathbf{g}\in \R^{n{\times}n}}\|\mathbf{m} - \mathbf{p}_{-s} \star \mathbf{g}\|^2_F + \alpha \|\mathbf{g}\|^2_F
\end{eqnarray}
where $\alpha>0$ is the regularization parameter.

We compute reconstructions using a FFT-based algorithm \cite{hansen2006deblurring}.

\subsection{TV regularization}

  We use the {\em isotropic} definition of total variation so the TV-regularized solution of \eqref{genrealregularrecon} and \eqref{genrealshiftedrecon} is defined on the two grids by:
\begin{eqnarray}\label{def:discreteTVregu}
{\mathbf{f}}^{(\alpha)} &=& \argmin_{\mathbf{f}\in \R^{n{\times}n}}\|\mathbf{m} - \mathbf{p}\star \mathbf{f}\|^2_F + \alpha \sum_{i,j=1}^n\sqrt{(D_x \mathbf{f})_{i,j}^2 + (D_y \mathbf{f})_{i,j}^2}, \\
\label{def:discreteTVshift}
{\mathbf{g}}^{(\alpha)} &=& \argmin_{\mathbf{g}\in \R^{n{\times}n}}\|\mathbf{m} - \mathbf{p}_{-s} \star \mathbf{g}\|^2_F + \alpha \sum_{i,j=1}^n \sqrt{(D_x \mathbf{g})_{i,j}^2 + (D_y \mathbf{g})_{i,j}^2}
\end{eqnarray}
where $\alpha>0$ is the regularization parameter and $D_x$, $D_y$ denote the horizontal and vertical partial derivative operators, respectively. 

We compute TV reconstructions using a standard primal-dual deconvolution solver \cite{chambolle2011first}. 
 

\subsection{Simulated Data (image with various geometric shapes)}\label{simudata}

\begin{figure}[H]
    \centering
   \includegraphics[width=350pt]{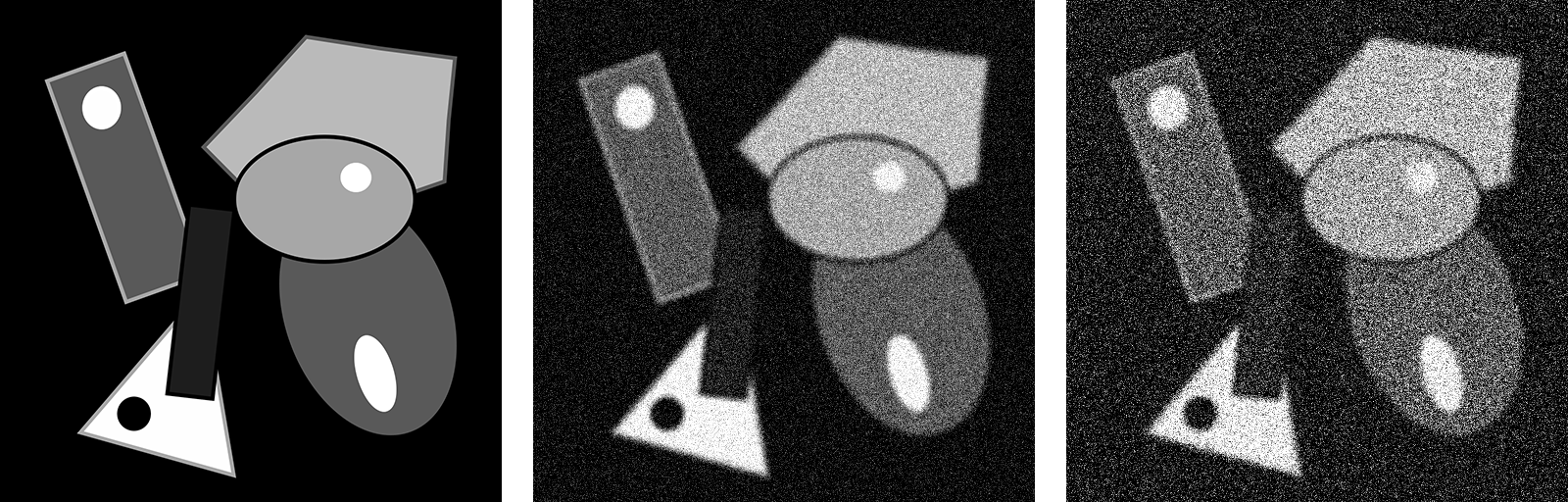}
    \caption{Left: Simulated ground truth with no noise or blur. Middle: Blurred with radius 4 kernel and 4\% added noise. Right: Blurred with radius 4 kernel and added 8\% noise.}
    \label{fig:manyshapes_gt}
\end{figure}

Our test image is a 512x512 sized computer created image with various geometric shapes with edges in several directions and various gray scale values for the shapes. See Figure \ref{fig:manyshapes_gt}.

The simulated image blur was created using a convolution kernel that has a constant value inside a disc domain and zero outside. We added Gaussian noise that is scaled relative to the norm of the image. Two different noise levels were chosen, one with 4\% and another with 8\% of noise relative to the image.

\subsection{Real Photographic Data}

To test our method on real data, we captured two small datasets with a Canon EOS 5D Mark IV camera. For the playing card set a Canon EF 100mm f/2.8 USM Macro lens was attached to the camera, and for the books dataset a newer Canon EF 100mm f/2.8 USM IS Macro lens was used. For both datasets multiple blur (lens misfocus) levels and noise (ISO) levels were captured. The scene was properly exposed each time. To compensate for the change in signal intensity/exposure when increasing the ISO setting the shutter speed was changed accordingly. The camera was mounted on a tripod and the subjects static so the change of shutter speed did not introduce any additional motion blur. For the playing card image the ISO levels were ISO100 for the "ground truth", ISO1600 for the less noisy image and ISO6400 for the more noisy image. For the image with books the ISO levels were ISO100, ISO6400 and ISO25600 accordingly. In both sets a blurred and noiseless (ISO100) image was used to estimate the amount of blur i.e. to estimate the radius of the blur kernel. In section \ref{dual-grid_tests} we will show tests of our method on the following real data examples.

\subsection{Real Dataset 1: Playing card}\label{queendata}

\begin{figure}[H]
    \centering
   \includegraphics[width= 340pt]{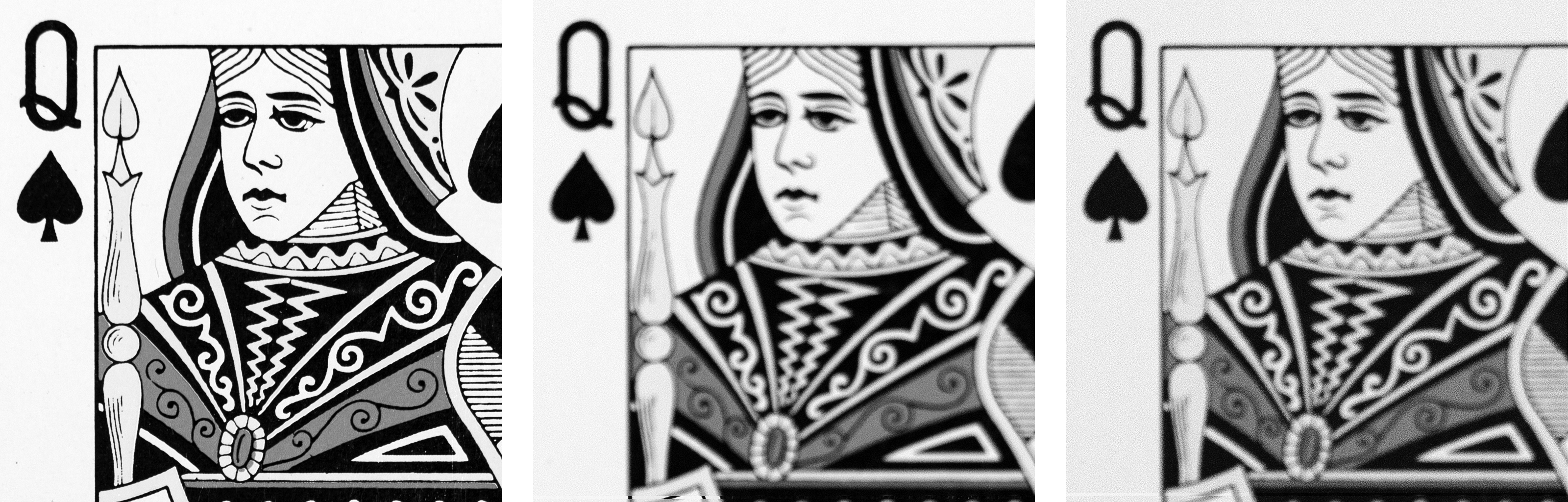}
    \caption{Left: "Ground truth", almost noiseless (ISO100). Middle: Slightly blurred, some noise (ISO1600). Right: Slightly blurred, more noise (ISO6400).}
    \label{fig:queentriplet}
\end{figure}

A 1024x1024 pixel sized crop from an image of a regular playing card. We have a ground truth estimate on the left in figure \ref{fig:queentriplet}. This is a real photograph correctly focused and taken with the lowest possible ISO setting. In the middle we have a slightly blurred and noisier version (ISO1600) of the same scene and on the right the same amount of blur but higher noise level (ISO6400).

\subsection{Real Dataset 2: Books}\label{booksdata}

\begin{figure}[H]
    \centering
   \includegraphics[width= 350pt]{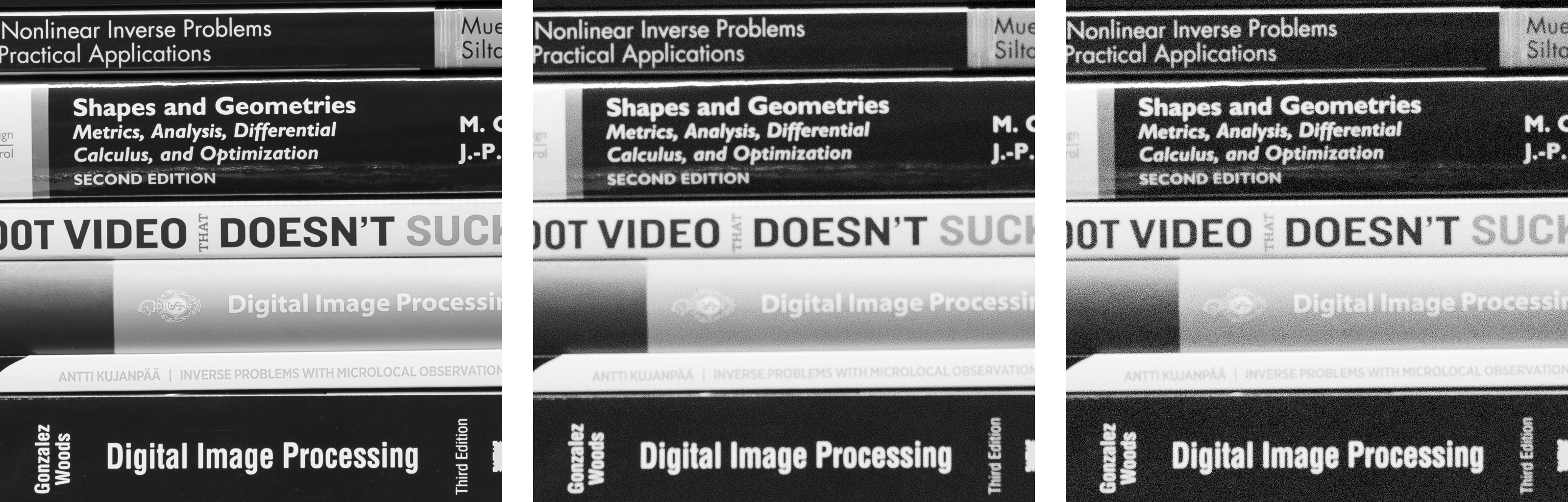}
    \caption{"Ground truth" (ISO100) on the left. Slightly blurred and small noise (ISO6400) in the middle. Slightly  blurred and more noise (ISO25600) on the right.}
    \label{fig:books_all}
\end{figure}

A 1024x1024 pixel sized crop from an image of a stack of books. We have a ground truth estimate on the left in Figure \ref{fig:books_all}. This is a real photograph correctly focused and taken with the lowest possible ISO setting. In the middle we have a slightly blurred and noisier version (ISO6400) of the same scene and on the right the same amount of blur but higher noise level (ISO25600).

\subsection{Natural image set} \label{naturalimageset}

\begin{figure}[H]
    \centering
   \includegraphics[width=350pt]{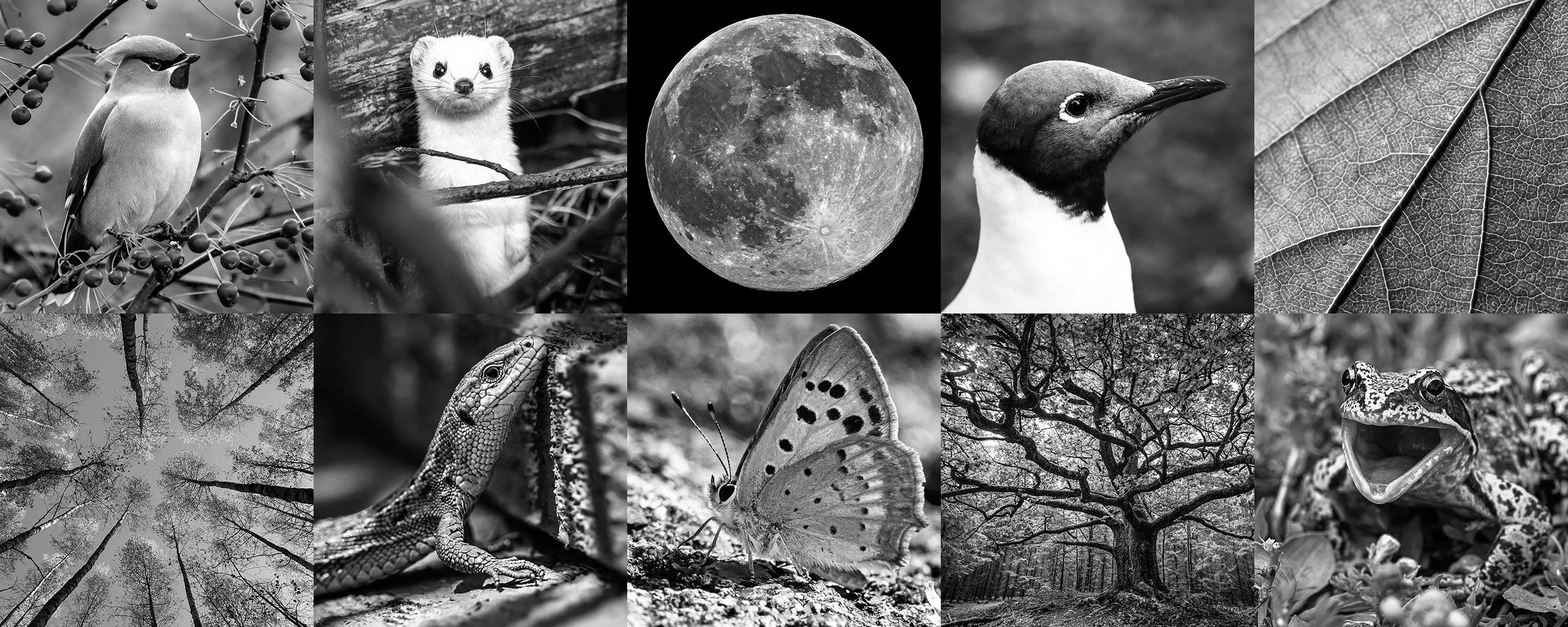}
    \caption{10 of the 230 natural images used for the numerical test of the similarity hypothesis. }
    \label{fig:collage_10of230}
\end{figure}

A set of 230 natural images photographed by the corresponding author during the last 10 years was used to carry out a few additional numerical experiments. These images have previously gone through various post-processing steps and were cropped or resized to be $1024\times 1024$ pixels for the purpose of this paper. Figure \ref{fig:collage_10of230} shows 10 examples of these 230 natural images. 

\section{Results}

\subsection{Numerical test of similarity hypothesis} \label{sec:hypo_backup}

\begin{figure}[H]
    \centering
   \includegraphics[width=350pt]{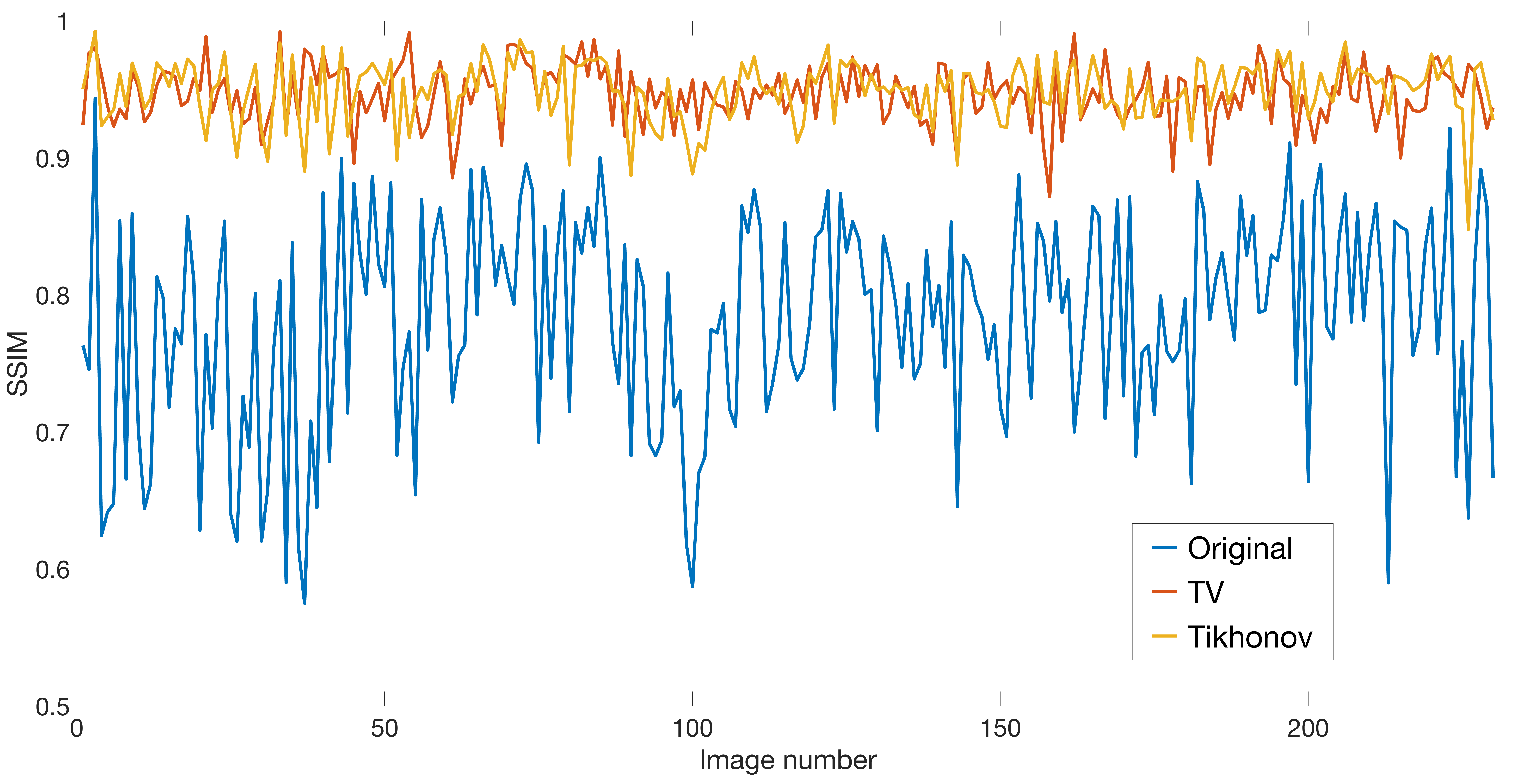}
    \caption{SSIM between images on the two grids. The blue graph indicates that for a regular photograph, a sub-pixel shift causes a significant change in the image as measured by SSIM. However, for images smoothed by Tikhonov or total variation denoising, the change caused by the shift is much weaker: their SSIM values are closer to 1. }
    \label{fig:ssim_graph}
\end{figure}

We first computed the SSIM between a set of 230 natural images $\mathbf{f}$ and their shifted versions $\mathbf{f}\star \mathbf{d}_s$. Figure \ref{fig:collage_10of230} shows 10 examples of the 230 natural image set. The blue graph in Figure \ref{fig:ssim_graph} shows that $\mathbf{f}$ is quite different from $\mathbf{f}\star \mathbf{d}_s$ under the SSIM measure. Then we repeated the calculation for denoised versions of those 230 images. (Note that we did not add any noise to the original images before denoising.) The red and orange graphs in Figure \ref{fig:ssim_graph} show that shifting an image smoothed by Tikhonov or TV regularization does not change the image as much in terms of SSIM.

\begin{figure}[H]
    \centering
   \includegraphics[width=350pt]{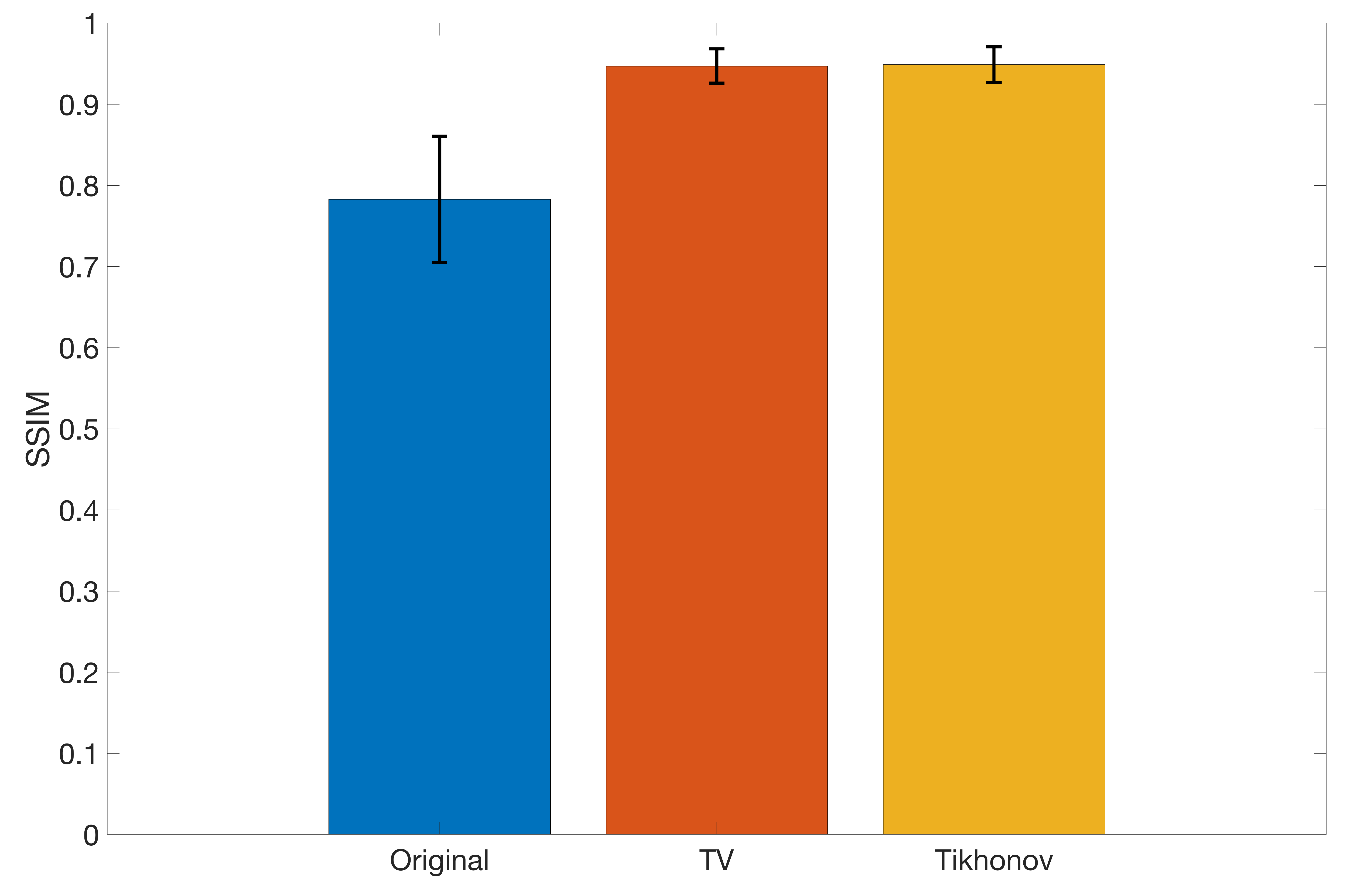}
    \caption{Mean SSIM value between images on the two grids and standard deviation computed for all 230 images in Figure \ref{fig:ssim_graph}.}
    \label{fig:ssim_bar}
\end{figure}

\clearpage 









\subsection{Dual-grid parameter selection tests}\label{dual-grid_tests}

In this section we present some results of the dual-grid parameter selection method for one simulated and two real photographic data sets. The simulated image data used for this test is presented in section \ref{simudata}. The real data sets are presented in sections \ref{queendata} and \ref{booksdata}.


\subsubsection{Simulated image example}

\begin{figure}[H]
    \centering
   \includegraphics[width=350pt]{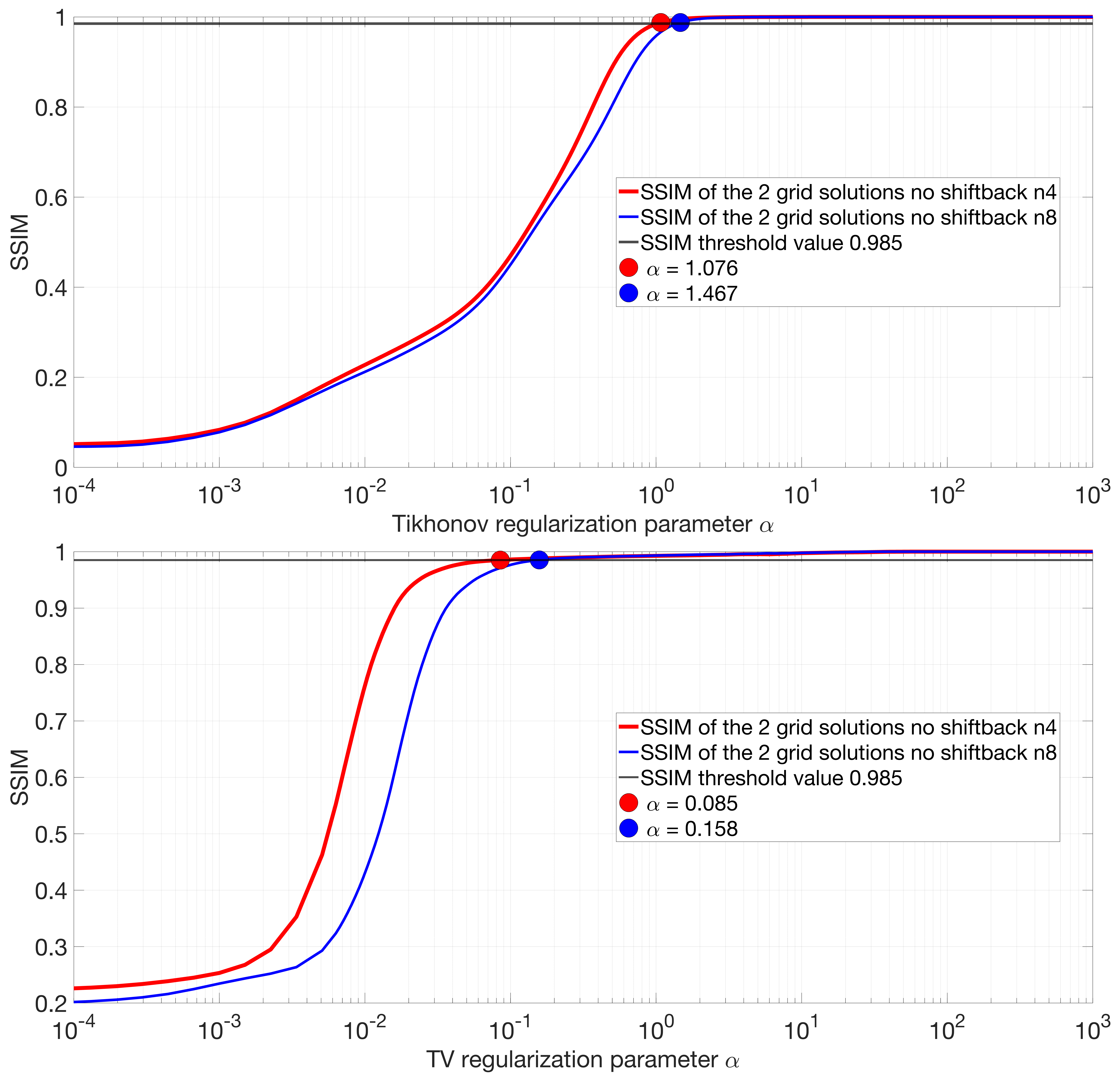}
    \caption{Dual-grid method for parameter choice for the case of simulated images shown in Figure \ref{fig:manyshapes_gt}. Plotted is the SSIM function between the two solutions SSIM$({\mathbf{g}}^{(\alpha)} ,{\mathbf{f}}^{(\alpha)})$ for two regularizers and two noise amplitudes. The SSIM threshold value is 0.985 in both cases. Top: Tikhonov regularization. The first regularization parameter after the selected threshold value is $\alpha = 1.076 $ for the case with less noise (red curve) and $\alpha = 1.467$ for the higher noise case (blue curve). Bottom: Total variation regularization. The first regularization parameter after the selected threshold value is $\alpha = 0.085 $ for the case with less noise (red curve) and $\alpha = 0.158$ for the higher noise case (blue curve).}
    \label{fig:manyshapesr4n4n8_ssimgraphs}
\end{figure}

First we test the dual-grid approach on a the simulated image example shown in Figure \ref{fig:manyshapes_gt}. For the simulated images we know that the blur kernel radius is 4 but the added noise makes solving this problem ill-posed. For both regularizers we get the expected behaviour of the SSIM function as seen in Figure \ref{fig:manyshapesr4n4n8_ssimgraphs}. We also see that for the lower noise case the selected parameter is in both cases smaller than with higher noise level.

\begin{figure}[H]
\begin{picture}(300,200)   
\put(0,0){\includegraphics[width=350pt]{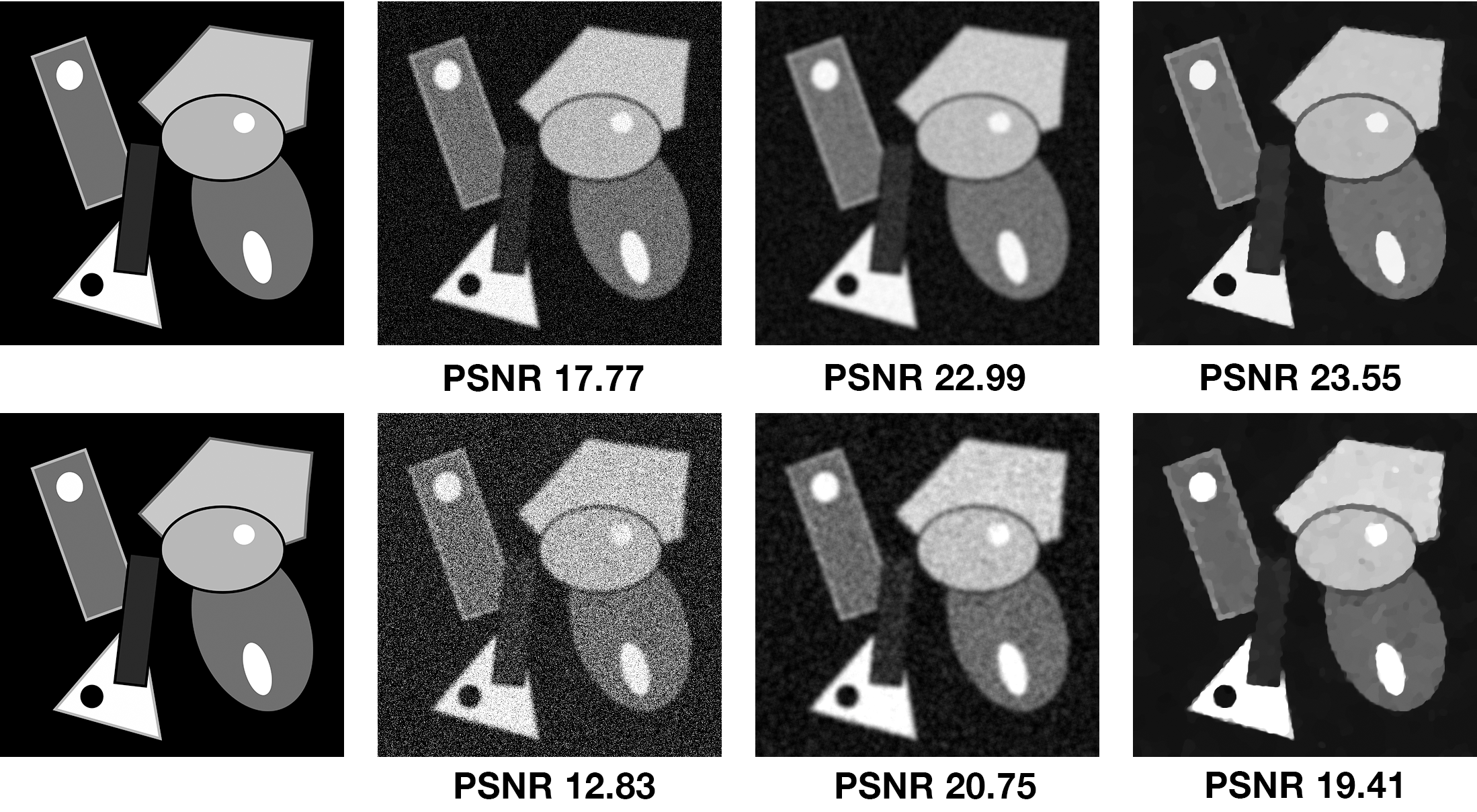}}
\put(5,200){Ground truth}
\put(92,200){Blurred and noisy}
\put(190,200){Tikhonov reg.}
\put(270,200){Total variation reg.}

    \end{picture}
    \caption{Geometric shapes images and reconstructions with regularization parameter values selected with dual-grid method as seen in Figure \ref{fig:manyshapesr4n4n8_ssimgraphs} . First row: Ground truth, blur radius 4 + noise 4\%, Tikhonov regularized, TV. Second row: Ground truth, blur radius 4 + noise 8\%, Tikhonov regularized , TV regularized (Tikh threshold .985, TV threshold .985). PSNR values between the ground truth and the noisy images and reconstructions are visible under each image.}
    \label{fig:shapes_r4n4_r4n8_tv_tikh}
\end{figure}

\subsubsection{Playing card example}

\begin{figure}[H]
    \centering
   \includegraphics[width=350pt]{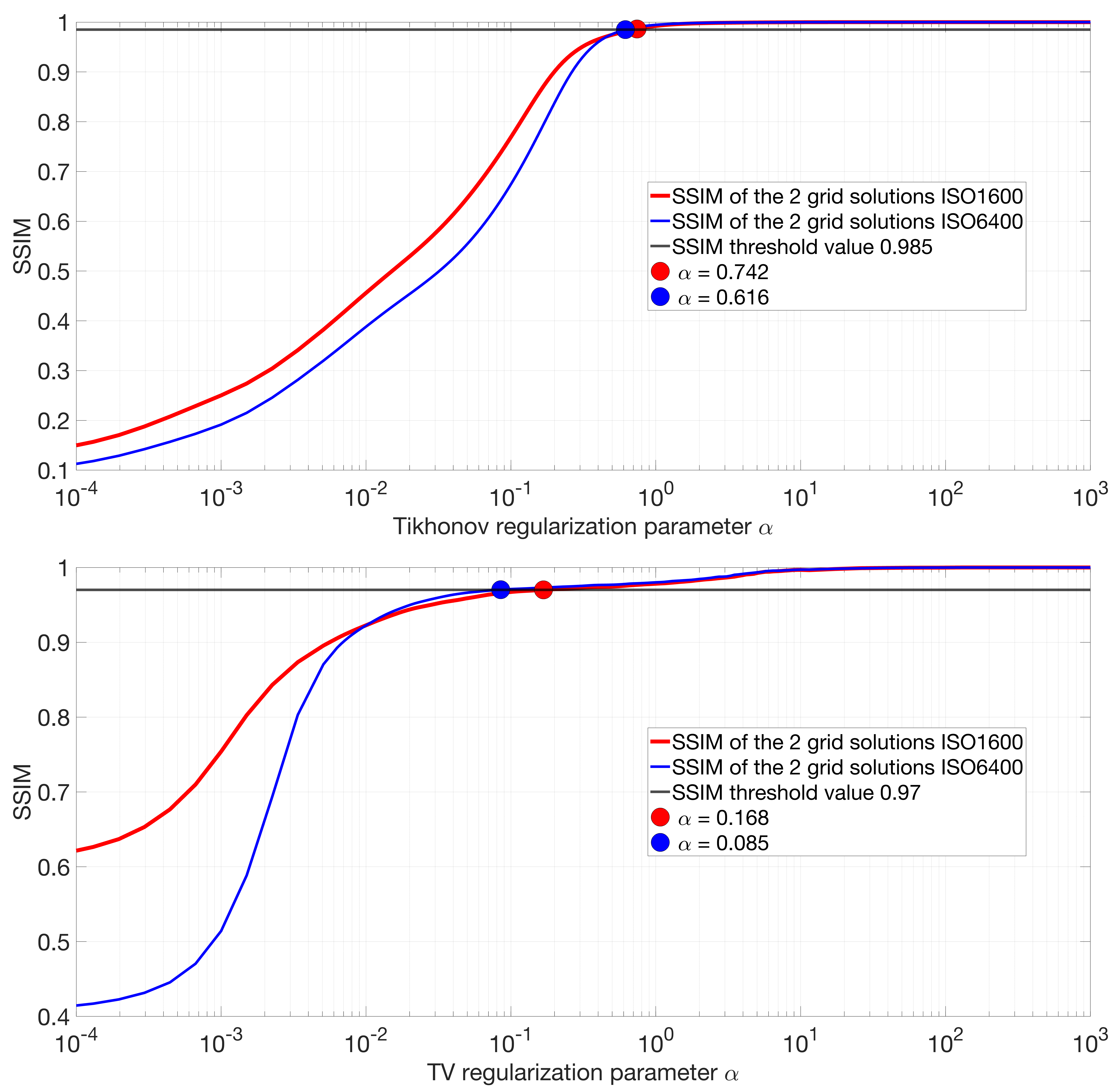}
    \caption{Dual-grid method for parameter choice for the case of the playing card images shown in Figure \ref{fig:queentriplet}. Plotted is the SSIM function between the two solutions SSIM$({\mathbf{g}}^{(\alpha)} ,{\mathbf{f}}^{(\alpha)})$ for two regularizers and two noise amplitudes. Top: Tikhonov regularization. The first regularization parameter after the selected SSIM threshold value of 0.985 is $\alpha = 0.742 $ for the case with less noise (red curve) and $\alpha = 0.616$ for the higher noise case (blue curve). Bottom: Total variation regularization. The first regularization parameter after the selected SSIM threshold value of 0.97 is $\alpha = 0.168 $ for the case with less noise (red curve) and $\alpha = 0.085$ for the higher noise case (blue curve).}
    \label{fig:queens_tikh_TV_ssimgraphs}
\end{figure}

For the first real data test, we chose the playing card image shown in Figure \ref{fig:queentriplet}. For the real data cases, we don't know the exact blur kernel that has caused the blur, so we need to estimate it. In this case, we ended up with an estimated kernel radius of 6 by measuring how much distinct points have spread after the misfocusing compared to the sharp noisless "ground truth" image. For the Tikhonov regularizer, we see pretty nicely behaving SSIM functions in Figure \ref{fig:queens_tikh_TV_ssimgraphs}. For the TV regularized ones, we see that the exact point where the SSIM function starts to flatten out is not so well defined. For both regularizers, we notice that we actually get a smaller parameter value choice for the higher noise case.

\begin{figure}[H]
\begin{picture}(320,180)   
\put(0,0){\includegraphics[width=350pt]{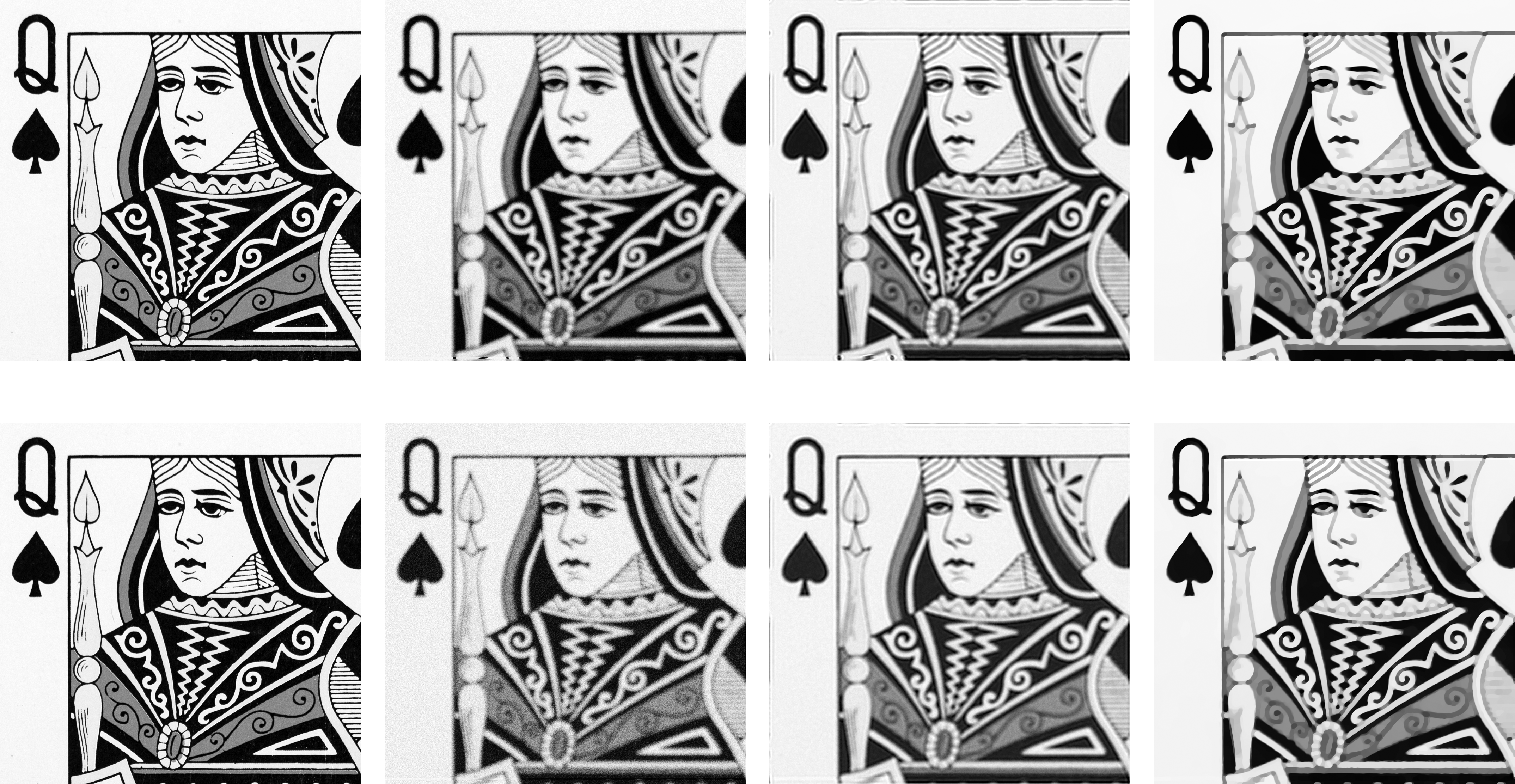}}
\put(10,185){Sharp image}
\put(92,185){Blurred and noisy}
\put(190,185){Tikhonov reg.}
\put(270,185){Total variation reg.}

    \end{picture}
    \caption{Queen playing card and reconstructions with regularization parameter values selected with dual-grid method as seen in Figure \ref{fig:queens_tikh_TV_ssimgraphs}. First row: "Ground truth" sharp image (ISO100), blur+little noise(ISO1600),  Tikhonov, TV. Second row: "Ground truth" sharp image (ISO100), blur+more noise(ISO6400), Tikhonov, TV. (Tikh threshold .985, TV threshold .97)}
    \label{fig:Queens_1600_6400_tv_tikh}
\end{figure}

Figure \ref{fig:Queens_1600_6400_tv_tikh} shows reconstructions of the blurry and noisy playing card images for the regularization parameter value estimations using the dual-grid approach for the selected thresholds as seen in Figure \ref{fig:queens_tikh_TV_ssimgraphs} (0.985 for Tikh. and 0.97 for TV).

\begin{figure}[H]
\begin{picture}(320,180)   
\put(0,0){\includegraphics[width=350pt]{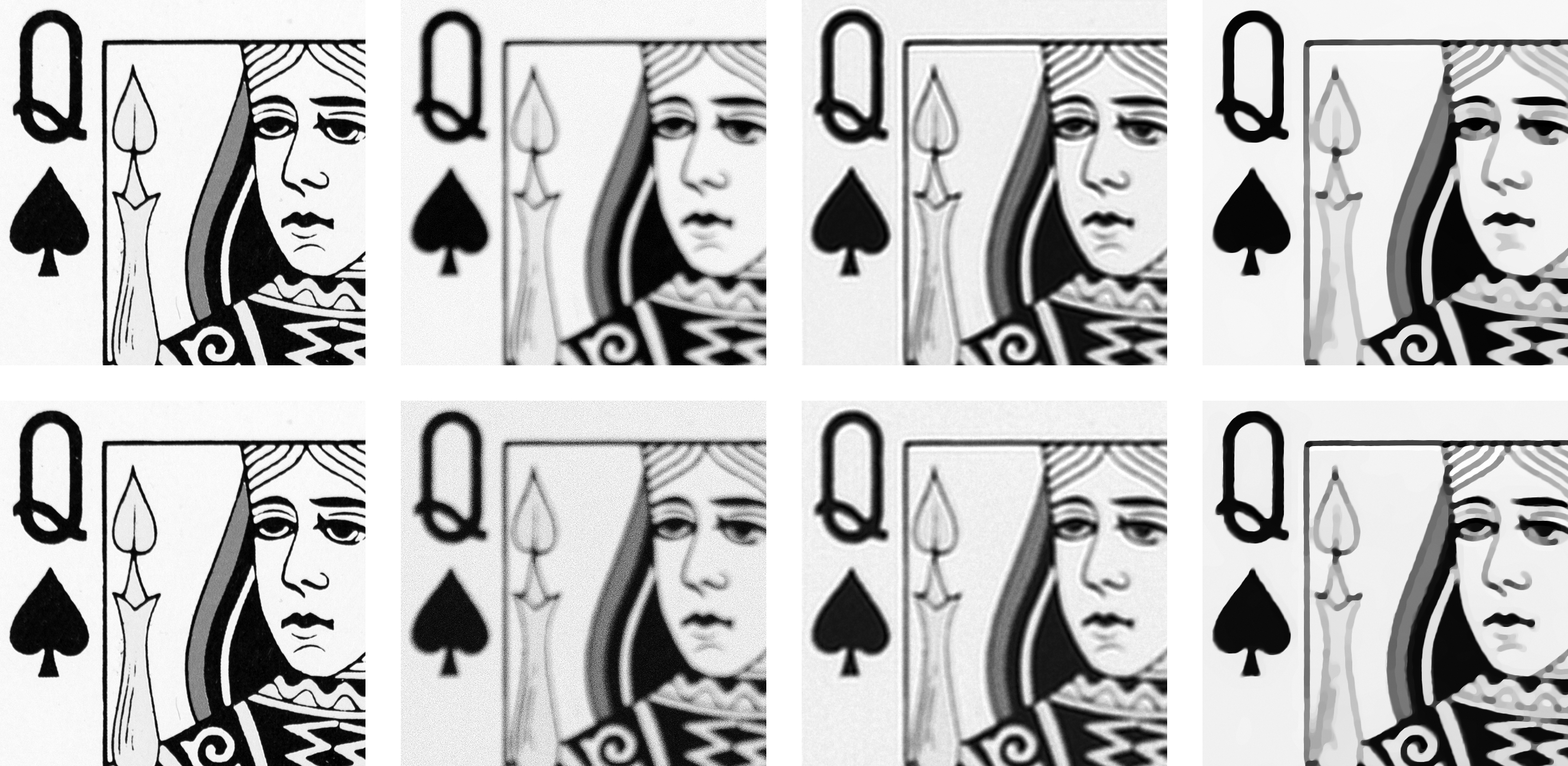}}
\put(10,185){Sharp image}
\put(92,185){Blurred and noisy}
\put(190,185){Tikhonov reg.}
\put(270,185){Total variation reg.}

    \end{picture}
    \caption{Queen playing card images from Figure \ref{fig:Queens_1600_6400_tv_tikh} cropped for a more detailed view. First row: "Ground truth" sharp image (ISO100), blur+little noise(ISO1600),  Tikhonov, TV. Second row: "Ground truth" sharp image (ISO100), blur+more noise(ISO6400), Tikhonov, TV. (Tikh threshold .985, TV threshold .97)}
    \label{fig:Queens_1600_6400_tv_tikh_crop}
\end{figure}

Figure \ref{fig:Queens_1600_6400_tv_tikh_crop} shows a cropped region of the larger images for a closer look at the details of the reconstructions.

\clearpage

\subsubsection{Books image example}

\begin{figure}[H]
    \centering
   \includegraphics[width=350pt]{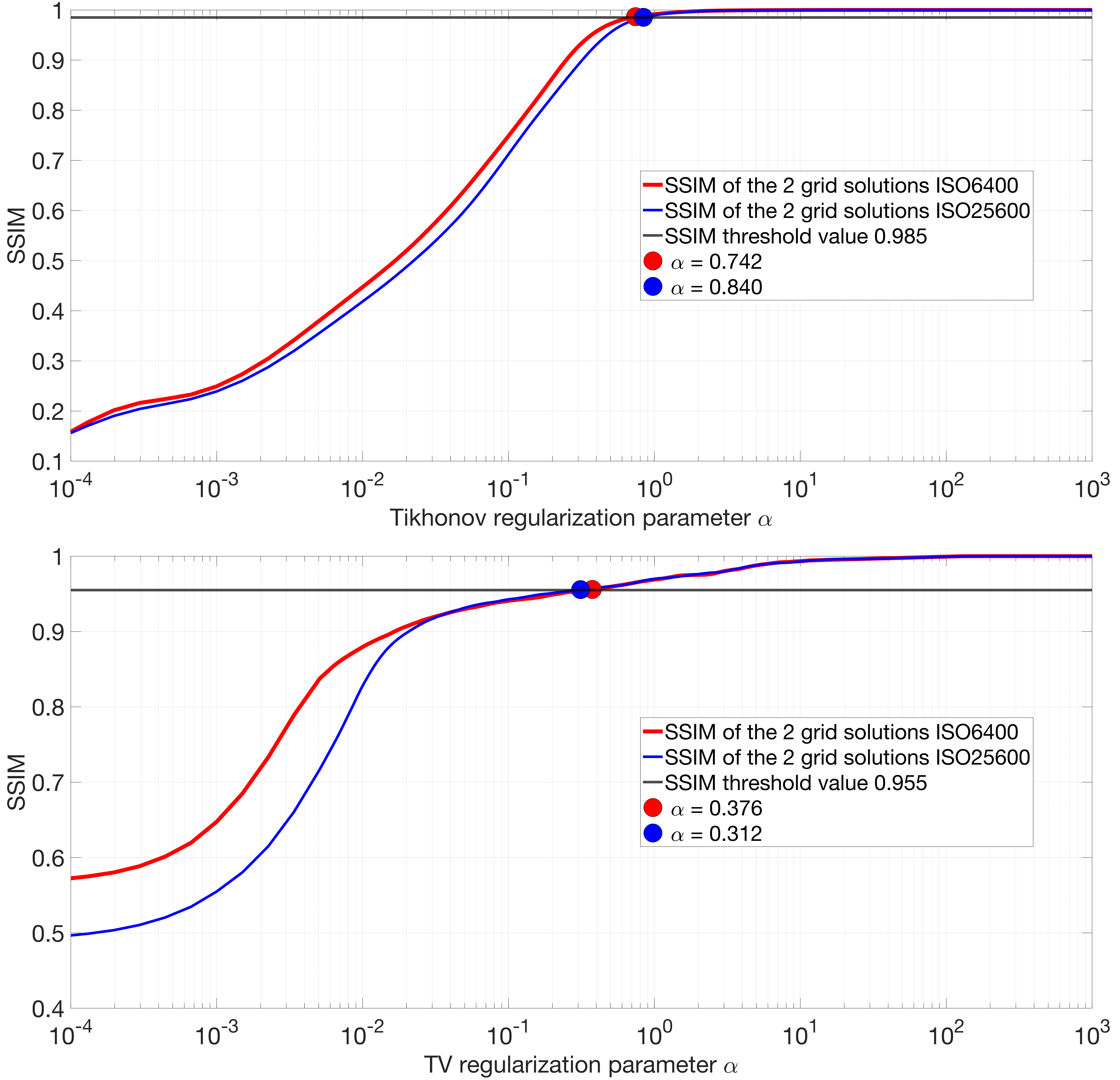}
    \caption{Dual-grid method for parameter choice for the case of the books images shown in Figure \ref{fig:books_all}. Plotted is the SSIM function between the two solutions SSIM$({\mathbf{g}}^{(\alpha)} ,{\mathbf{f}}^{(\alpha)})$ for two regularizers and two noise amplitudes. Top: Tikhonov regularization. The first regularization parameter after the selected threshold value of 0.985 is $\alpha = 0.742 $ for the case with less noise (red curve) and $\alpha = 0.84$ for the higher noise case (blue curve). Bottom: Total variation regularization. The first regularization parameter after the selected threshold value of 0.955 is $\alpha = 0.376 $ for the case with less noise (red curve) and $\alpha = 0.312$ for the higher noise case (blue curve).}
    \label{fig:books_tikh_TV_ssimgraphs}
\end{figure}

As the second case of real data we have the Books images seen in Figure \ref{fig:books_all}. The blur kernel radius estimation in this case is 5. We see a similar behavior of the SSIM functions to the playing card image case.
In the TV case, we clearly notice how the SSIM log plots in Figure \ref{fig:books_tikh_TV_ssimgraphs} don't have the same s-shape, but they grow more gradually from a value of 0.9 to 1. This makes selecting the threshold more challenging.

\begin{figure}[H]
\begin{picture}(300,180)   
\put(0,0){\includegraphics[width=350pt]{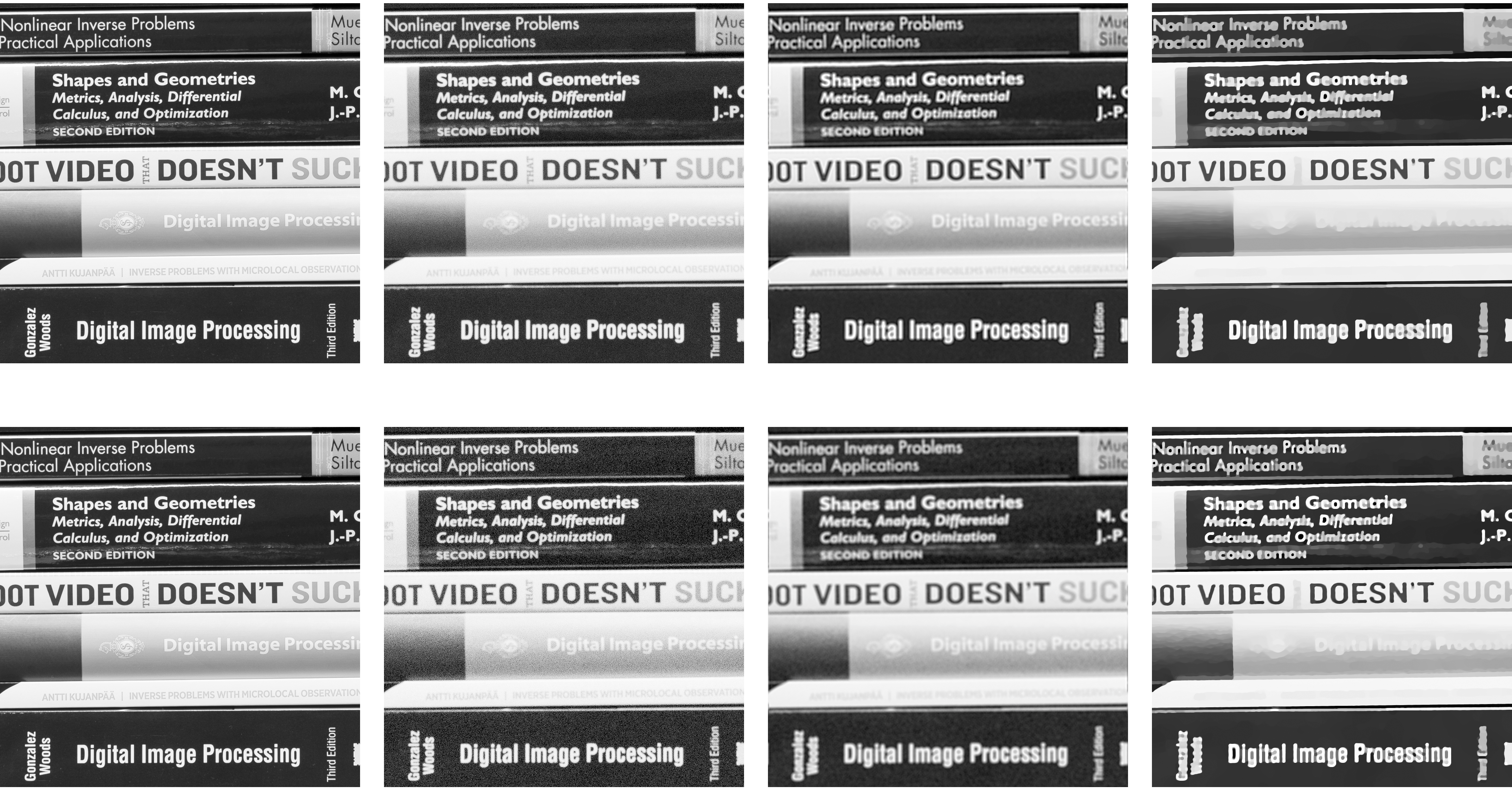}}
\put(5,188){Sharp image}
\put(92,188){Blurred and noisy}
\put(190,188){Tikhonov reg.}
\put(268,188){Total variation reg.}

    \end{picture}
    \caption{Books images and reconstructions with regularization parameter values selected with dual-grid method as seen in Figure \ref{fig:books_tikh_TV_ssimgraphs}. First row: "Ground truth" sharp image (ISO100), Blur+little noise, Tikhonov, TV. Second row: "Ground truth" sharp image (ISO100), blur+more noise, Tikhonov, TV (Tikh threshold .985, TV threshold .955)}
    \label{fig:Books_6400_25600_tikh_tv}
\end{figure}

\begin{figure}[H]
\begin{picture}(300,180)   
\put(0,0){\includegraphics[width=350pt]{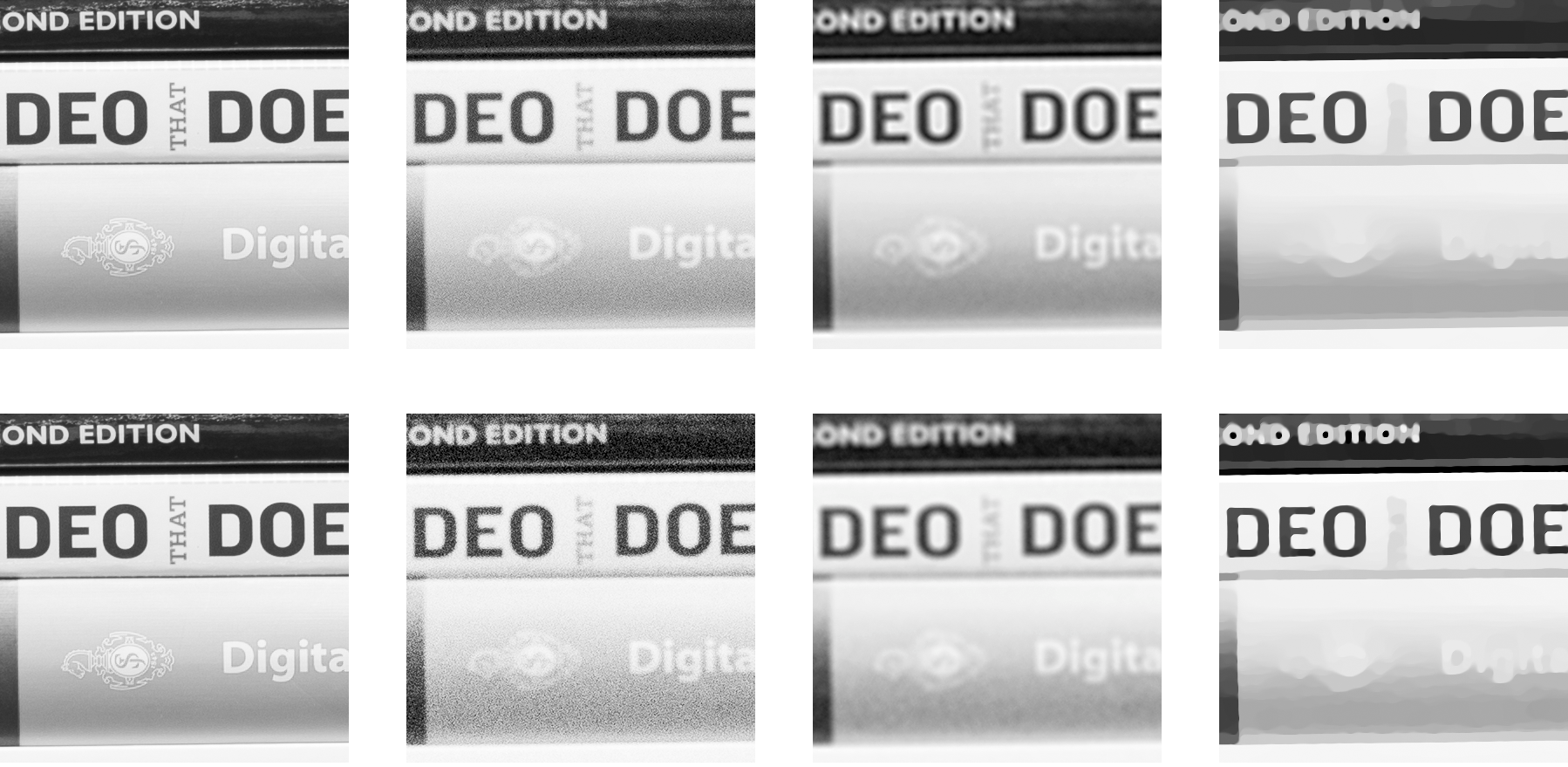}}
\put(5,180){Sharp image}
\put(92,180){Blurred and noisy}
\put(190,180){Tikhonov reg.}
\put(270,180){Total variation reg.}

    \end{picture}
    \caption{Books images from Figure \ref{fig:Books_6400_25600_tikh_tv} cropped for more detailed look. First row: "Ground truth" sharp image (ISO100), blur+little noise, Tikhonov, TV. Second row: "Ground truth" sharp image (ISO100), blur+more noise, Tikhonov, TV (Tikh. threshold .985, TV threshold .955)}
    \label{fig:Books_6400_25600_tv_tikh_zoom}
\end{figure}

Figure \ref{fig:Books_6400_25600_tikh_tv} shows reconstructions of the blurry and noisy Book images for the regularization parameter value estimations using the dual-grid approach for the selected thresholds as seen in Figure \ref{fig:books_tikh_TV_ssimgraphs} (0.985 for Tikh. and 0.955 for TV).
Figure \ref{fig:Books_6400_25600_tv_tikh_zoom} shows a cropped region of the images for a closer look at the details of the reconstructions.

\subsection{Comparison with discrepancy principle and bilevel optimization parameter choice}\label{discrep_comparison}

To have a comparison to other parameter choice methods we first use the discrepancy principle for Tikhonov and TV, and second a bilevel optimization method in the TV case.

The goal in using the discrepancy principle is to find the parameter value $\alpha$ so that the regularized solution $\mathbf{f}^{(\alpha)}$ has a residual norm comparable to the noise level $\delta$ in the data $\mathbf{m}$. For this we define 
\begin{equation}\label{discr_Psi}
    \Psi(\alpha) := \|\mathbf{p}\star \mathbf{f}^{(\alpha)} - \mathbf{m}\|_F
\end{equation} 
and take the discrepancy-based parameter to be the $\alpha$ for which $\Psi(\alpha)\approx \delta $.

 For the noise estimation we first sampled 4 constant valued areas of each of the test images and computed the standard deviation $\sigma$ of the pixel values in these regions. Then we calculated the mean of these values. To get the final noise estimate $\delta$ we multiply the mean std value with the square root of the number of pixels in that image, $\delta = \sigma_{mean} \sqrt{n\times n} $. We can then plot the residual norm for each $\alpha$ and see when it matches the noise estimate.

\begin{figure}[H]
    \centering
   \includegraphics[width=350pt]{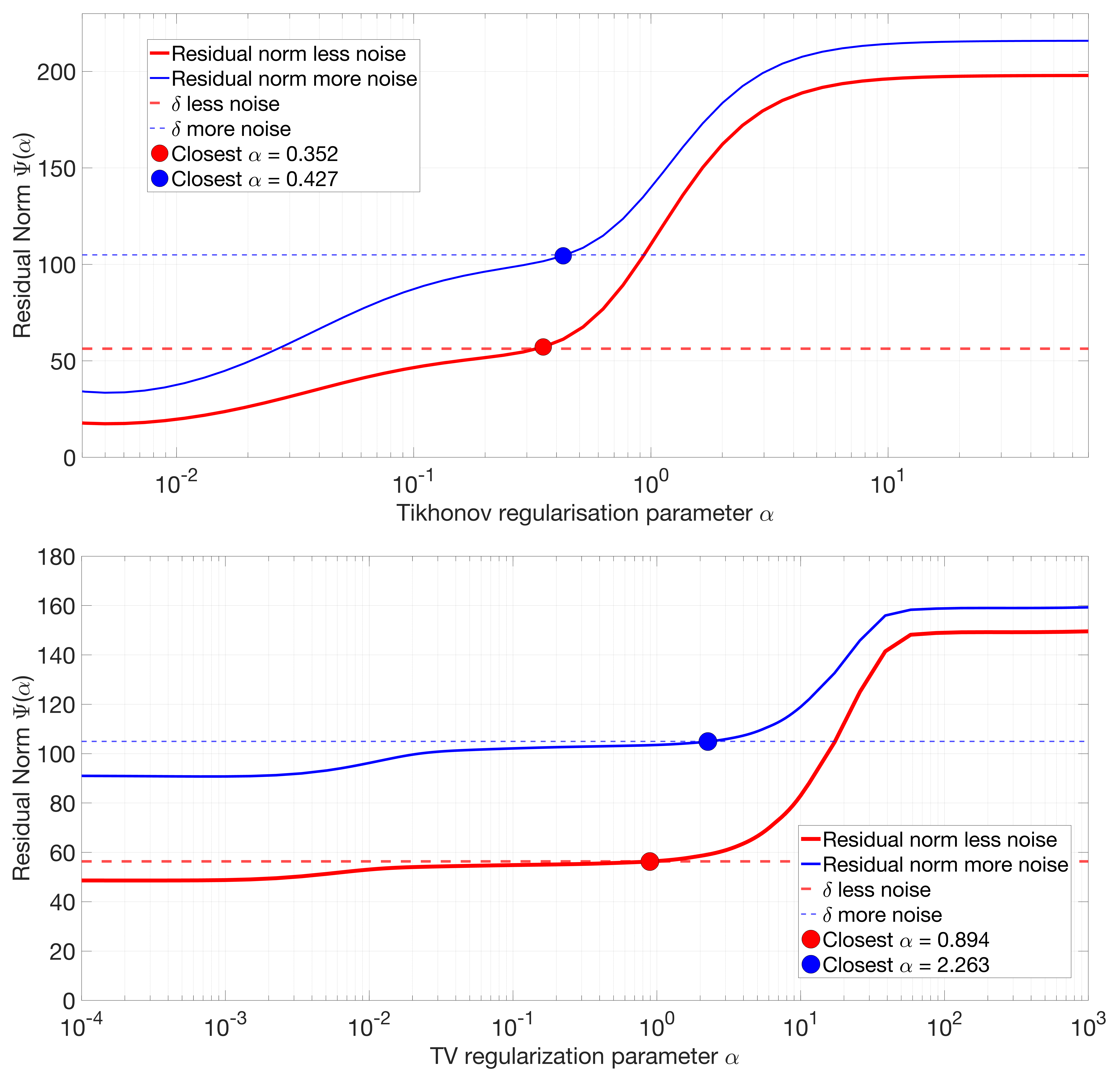}
    \caption{Discrepancy principle for parameter choice for the case of simulated images shown in Figure \ref{fig:manyshapes_gt}. Plotted is the function $\Psi(\alpha)$ defined in (\ref{discr_Psi}) for two regularizers and two noise amplitudes. Top: Tikhonov regularization. The closest regularization parameter to the estimated noise level $\delta$ is $\alpha = 0.352 $ for the case with less noise (red curve) and $\alpha = 0.427$ for the higher noise case (blue curve). Bottom: Total variation regularization. The closest regularization parameter to the estimated noise level $\delta$ is $\alpha = 0.821 $ for the case with less noise (red curve) and $\alpha = 1.894$ for the higher noise case (blue curve).}
    \label{fig:shapes_Tikh_TV_discrep_plots}
\end{figure}

\begin{figure}[]
    \centering
   \includegraphics[width=300pt]{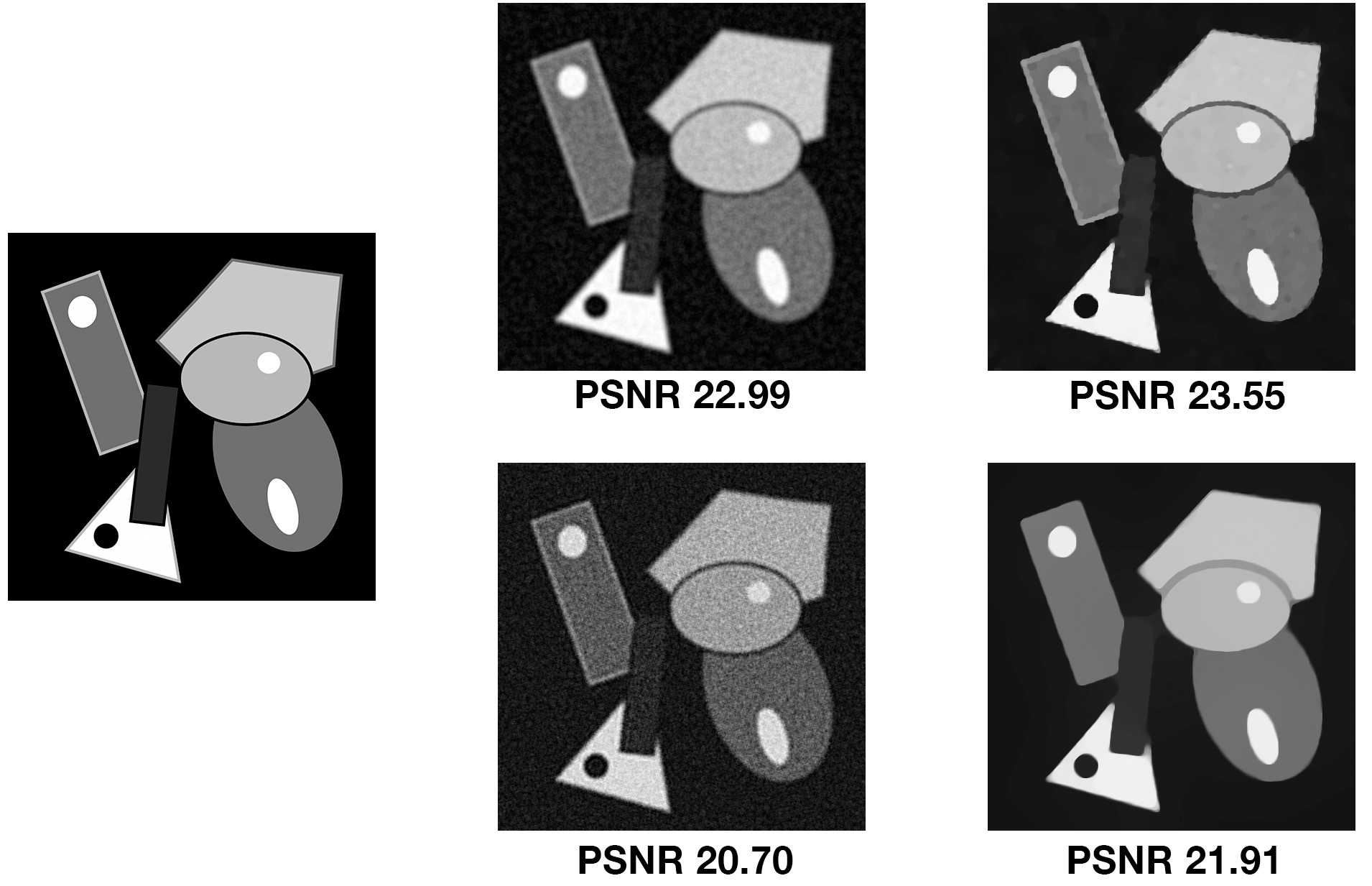}
    \caption{Comparison of optimal reconstructions according to dual-grid method and discrepancy principle for the simulated image with noise level 4\%. Left: Ground truth. Upper row: Dual-grid results. Bottom row: Discrepancy principle results. PSNR values compared to the ground truth image are under each reconstruction.}
    \label{fig:shapes_Tikh_TV_dual_discrep_comparison_n4}
\end{figure}

\begin{figure}[]
    \centering
   \includegraphics[width=300pt]{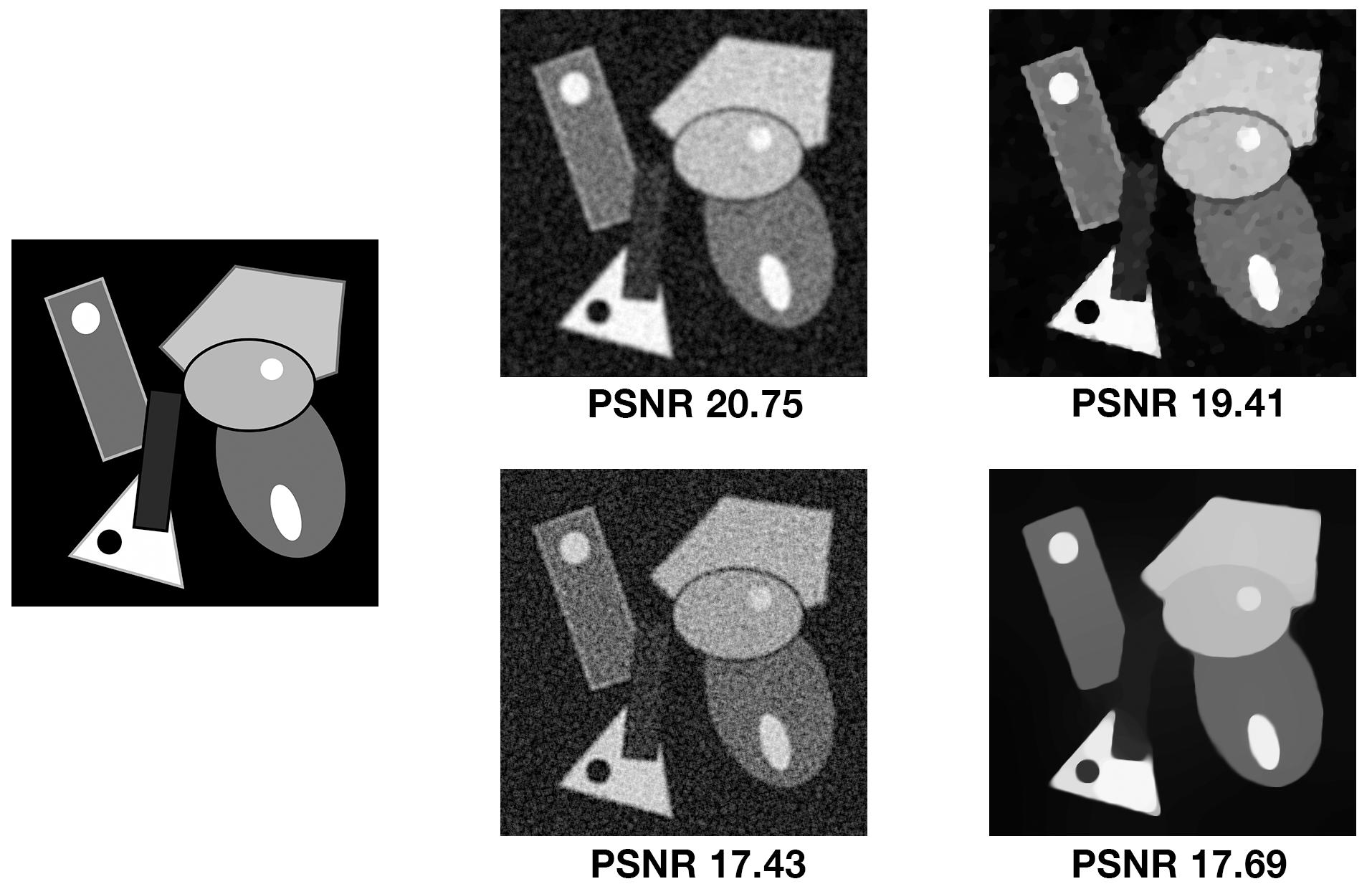}
    \caption{Comparison of optimal reconstructions according to dual-grid method and discrepancy principle for the simulated image with noise level 8\%. Left: Ground truth. Upper row dual-grid. Lower row discrepancy principle. PSNR values compared to the ground truth image are under each reconstruction.}
    \label{fig:shapes_Tikh_TV_dual_discrep_comparison_n8}
\end{figure}

The discrepancy plots for the simulated image example in Figure \ref{fig:shapes_Tikh_TV_discrep_plots} show clear intersections with the estimated noise for both noise levels in both the Tikhonov regularization case and the total variation regularization case. The according reconstructions compared to the dual-grid method can be seen in Figure \ref{fig:shapes_Tikh_TV_dual_discrep_comparison_n4} for the lower noise case and in Figure \ref{fig:shapes_Tikh_TV_dual_discrep_comparison_n8} for the higher noise case. In both Tikhonov regularization cases the discrepancy principle gives us a lower parameter value than the proposed dual-grid method and for Total variation regularization the discrepancy principle gives us higher parameter values in comparison to the dual-grid method. All the parameter values are also visible in Table \ref{alphas_tableTikh}.

\begin{figure}[H]
    \centering
   \includegraphics[width=340pt]{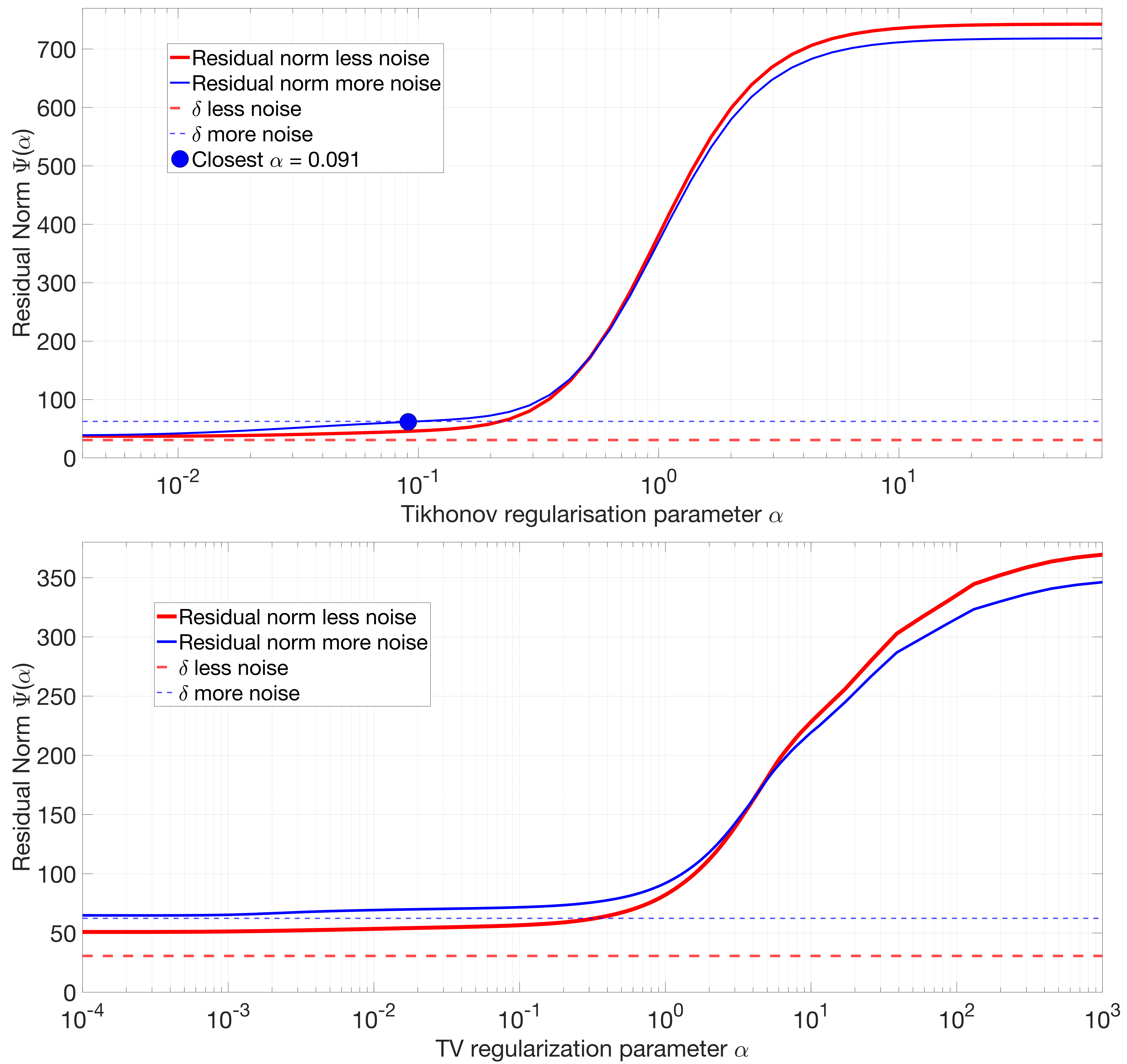}
    \caption{Discrepancy principle for parameter choice for the case of the playing card images shown in Figure \ref{fig:queentriplet}. Plotted is the function $\Psi(\alpha)$ defined in (\ref{discr_Psi}) for two regularizers and two noise amplitudes. Top: Tikhonov regularization. We only get a parameter value for the higher noise case. The closest value to the intersection with the estimated noise level is $\alpha = 0.091$ for the higher noise case (blue curve). Bottom: Total variation regularization. The discrepancy principle fails to give us parameter values for the estimated noise levels. }
    \label{fig:queen_both_Tikh_TV_discrep}
\end{figure}

The discrepancy plots for the playing card image example intersect with the estimated noise level only for the Tikhonov regularized case for the higher noise level in Figure \ref{fig:queen_both_Tikh_TV_discrep}. The lower noise Tikhonov curve and both total variation regularization cases fail to give us a parameter choice as seen in the same Figure \ref{fig:queen_both_Tikh_TV_discrep}. The parameter value is significantly smaller than the one we get with the dual-grid method. All parameter values can be compared in Table \ref{alphas_tableTikh}.


\begin{figure}[H]
    \centering
   \includegraphics[width=340pt]{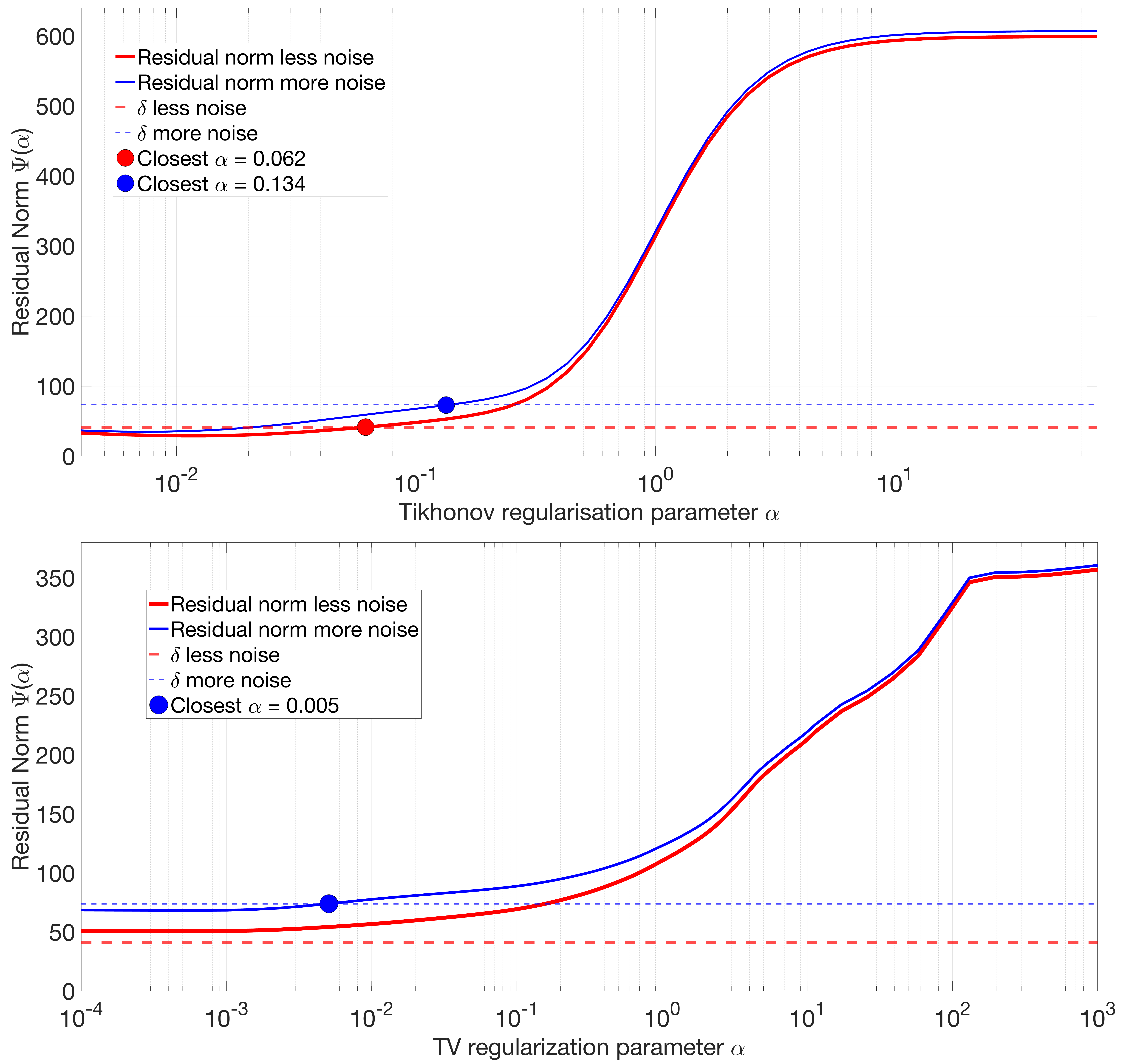}
    \caption{Discrepancy principle for parameter choice for the case of the books images shown in Figure \ref{fig:books_all}. Plotted is the function $\Psi(\alpha)$ defined in (\ref{discr_Psi}) for two regularizers and two noise amplitudes. Top: Tikhonov regularization. The closest regularization parameter to the estimated noise level $\delta$ is $\alpha = 0.062 $ for the case with less noise (red curve) and $\alpha = 0.134$ for the higher noise case (blue curve). Bottom: total variation regularization. We fail to find a parameter value for the case with less noise (red curve). The closest regularization parameter to the estimated noise level $\delta$ is $\alpha = 0.005$ for the higher noise case (blue curve).}
    \label{fig:books_tikh_TV_discrep}
\end{figure}

The discrepancy plots for the books image example in Figure \ref{fig:books_tikh_TV_discrep} show intersections with the estimated noise for both noise levels in the Tikhonov regularization case seen. In the total variation regularization case we get a $\alpha$ value but only in the larger noise case. All the parameter values can be seen in Table \ref{alphas_tableTikh}.

\begin{table}[!htbp]
\centering
\begin{tabular}{|c|c|c|c|c|c}
	\hline
	 Image test &  Tikh. discrep.$\alpha$ & \textbf{Tikh. Dual-grid $\alpha$}  \\
	\hline\hline
	Simulated 4\% & 0.352 & 1.076  \\
	\hline
	Simulated 8\% & 0.427  & 1.467  \\
	\hline
	Queen low     & -     &  0.742  \\
	\hline
	Queen high    & 0.091 & 0.616  \\
	\hline
	Books low     & 0.062 &  0.742    \\
	\hline
	Books high    & 0.134 &  0.840   \\
	\hline
\end{tabular}
\caption{Optimal Tikhonov regularization parameter values $\alpha$ for all image examples according to the discrepancy principle and the proposed dual-grid method.}
\label{alphas_tableTikh}
\end{table}

For the bilevel optimization comparison, we used the FIFB algorithm introduced in \cite{suonpera2024linearly} with implementations available on Zenodo \cite{suonpera_2023_7974062}. Several things to note about comparing our method with the bilevel optimization approach. Firstly, it needs a ground truth reference image, and we only have the sharp, low noise estimates for the real data available. Secondly, it uses smoothed TV, and lastly, it needs substantial computational resources to reach convergence. The resulting parameters in comparison with the proposed dual-grid method can be found in \ref{alphas_tableTV}. We see that for the simulated image, the alpha values are quite close to each other for both methods. For the playing card image, the bilevel alpha values are about half of the dual-grid ones. For the Books image the values for the lower noise case are quite close, but in the higher noise case the bilevel method gives a larger parameter this time. These inconsistencies with real data are most likely caused by the fact that we only have an estimate of the ground truth for the bilevel method and the SSIM threshold selection for the dual-grid method becomes difficult in some real data cases also.

\begin{table}[!htbp]
\centering
\begin{tabular}{|c|c|c|c|c}
	\hline
	 Image test & TV discrep.$\alpha$ & \textbf{TV Dual-grid $\alpha$} & TV Bilevel (FIFB) $\alpha$ \\
	\hline\hline
	Simulated 4\%  & 0.894 & 0.085 & 0.064  \\
	\hline
	Simulated 8\%  & 2.263 & 0.158 & 0.123 \\
	\hline
	Queen low      &  -  & 0.168 &  0.085 \\
	\hline
	Queen high     &  -  &  0.085 & 0.048  \\
	\hline
	Books low      & -  &  0.376   & 0.440  \\
	\hline
	Books high     & 0.005 & 0.312  & 0.559 \\
	\hline
\end{tabular}
\caption{Optimal TV parameter values $\alpha$ for all image examples. The discrepancy principle, the proposed dual-grid method and the FIFB bilevel optimization algorithm.}
\label{alphas_tableTV}
\end{table}

\begin{figure}[!htbp]
    \centering
    \includegraphics[width=350pt]{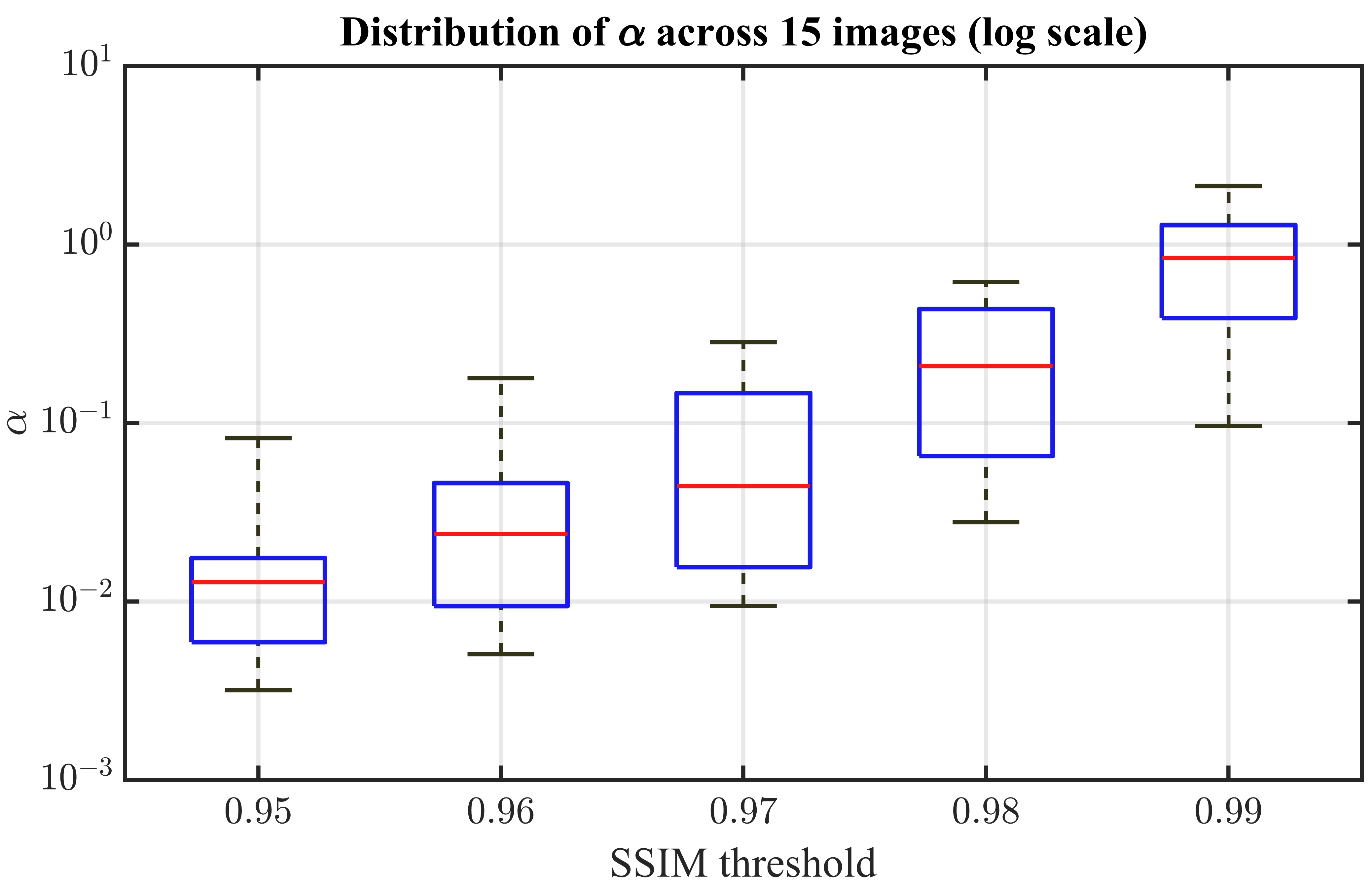}
    \caption{Box plots of $\alpha$ across 15 of the 230 natural images for each SSIM threshold
$\{0.95,0.96,0.97,0.98,0.99\}$. Boxes show the interquartile range (IQR = $[Q_1,Q_3]$),
the horizontal line is the median, whiskers extend to the most extreme values within
$1.5\times\mathrm{IQR}$.
The vertical axis is labeled as powers of ten, but boxes are computed on $\log_{10}(\alpha)$,
so equal vertical distances correspond to multiplicative factors in $\alpha$.
Medians increase monotonically with the threshold (0.0129, 0.0239, 0.0444, 0.2086, 0.8402).
The IQR also grows (0.0116, 0.0327, 0.1258,
0.3485, 0.8534), indicating substantial image-to-image variability at fixed thresholds.}
    \label{fig:15imboxplot}
\end{figure}

We did an additional experiment with 15 images of the 230 natural image set presented in Section \ref{naturalimageset}. First we blurred the images with a blurring kernel of radius 4 and added some noise (AWGN) with a standard deviation of 0.02. Then we used the dual-grid method to look for the regularization parameter with different SSIM thresholds from 0.95 to 0.99. The results are presented in Table \ref{tab:alpha_thresholds_simple_rounded}. 
To visualize how the selected regularization parameter varies across images and thresholds, we plot $\alpha$ on a logarithmic scale using box plots (Fig.~\ref{fig:15imboxplot}). Two consistent patterns emerge. First, $\alpha$ increases strictly with the SSIM threshold for every image, with the median rising from $0.0129$ at $0.95$ to $0.8402$ at $0.99$ (approximately $\times 65$). Second, within a fixed threshold the spread across images is large (IQRs $=\{0.0116,\,0.0327,\,0.1259,\,0.3485,\,0.8534\}$ for thresholds $\{0.95,\dots,0.99\}$), demonstrating that a single global $\alpha$ is not representative for a natural image set with different intensity level changes and various levels of details. Note that because the vertical scale is logarithmic, equal vertical distances correspond to multiplicative changes in $\alpha$, which is appropriate given the multiorder-of-magnitude range of values.

\begin{table}[!htbp]
\centering
\small
\setlength{\tabcolsep}{5pt}
\renewcommand{\arraystretch}{1.15}
\begin{tabular}{|l|c|c|c|c|c|}
\hline
 & \multicolumn{5}{c|}{\textbf{SSIM threshold}} \\
\cline{2-6}
\textbf{Image \#} & \textbf{0.95} & \textbf{0.96} & \textbf{0.97} & \textbf{0.98} & \textbf{0.99} \\
\hline\hline
Im1  & 0.00508 & 0.00808 & 0.015  & 0.0824 & 0.387 \\
\hline
Im2  & 0.00593 & 0.00944 & 0.015  & 0.0605 & 0.452 \\
\hline
Im3  & 0.015   & 0.0239  & 0.038  & 0.0824 & 0.209 \\
\hline
Im4  & 0.00593 & 0.00944 & 0.0175 & 0.0605 & 0.284 \\
\hline
Im5  & 0.0205  & 0.0518  & 0.153  & 0.452  & 1.34  \\
\hline
Im6  & 0.0326  & 0.0824  & 0.209  & 0.528  & 1.82  \\
\hline
Im7  & 0.0129  & 0.0326  & 0.131  & 0.528  & 2.13  \\
\hline
Im8  & 0.0110  & 0.0205  & 0.0444 & 0.244  & 1.15  \\
\hline
Im9  & 0.0175  & 0.0326  & 0.0706 & 0.209  & 0.720 \\
\hline
Im10 & 0.00319 & 0.00508 & 0.00944& 0.0279 & 0.0962 \\
\hline
Im11 & 0.00808 & 0.0150  & 0.0380 & 0.153  & 0.840 \\
\hline
Im12 & 0.0824  & 0.179   & 0.284  & 0.616  & 1.34  \\
\hline
Im13 & 0.00373 & 0.00593 & 0.0129 & 0.0518 & 0.387 \\
\hline
Im14 & 0.0129  & 0.0279  & 0.0824 & 0.284  & 0.981 \\
\hline
Im15 & 0.0175  & 0.0706  & 0.153  & 0.387  & 0.981 \\
\hline
\end{tabular}
\caption{TV deblurring parameter $\alpha$ per image for each SSIM threshold using the dual-grid method. Values rounded to three significant digits for readability.}
\label{tab:alpha_thresholds_simple_rounded}
\end{table}

\section{Discussion }

The parameter selection approach we propose is quite flexible in principle. The shifted grid idea can be used with any forward model, not only convolution.
The regularization term can also take many forms, including others than the Tikhonov and TV types considered here. The data fidelity norm and image similarity measures can be changed as well within the same parameter choice framework. Detailed assumptions on all of these choices while guaranteeing unique and useful parameter choices are outside the scope of this paper. 

A particular strength of our method is that the noise amplitude is not needed. This is a great advantage in real data cases. As we saw in the comparison section \ref{discrep_comparison}, estimating the noise amplitude may lead to difficulties with traditional parameter choice methods such as the discrepancy principle. 

Our experimental results show that the dual-grid method works for the image deblur problem with both simulated and real data. 
We do need to select the threshold value in our approach. This selection is influenced at least by the regularizer and the amount of fine textures in the image. With a fixed application and specific end-users' needs, a suitable threshold can be found by experimentation and then kept in subsequent work. 

The slight blurring effect of the shift operator does not cause problems as it is only used to determine an optimal parameter. In other words, one would not use the shifted reconstruction as the final solution as it has additional blur.

We noticed a curious effect in our computational experiments. For the real image data, the dual-grid method may give a larger regularization parameter for lower noise than for higher noise. This happened in both TV regularization experiments and for the playing card image with Tikhonov regularization. In the TV case $\alpha$ weights the sparsity of the gradient so the optimal $\alpha$ need not grow with the noise level and can remain roughly constant, or even decrease, because the prior term targets the same sparsity structure across noise levels. In this sense, the observation of smaller $\alpha$ at higher noise is not contradictory but consistent with TV’s sparsity prior. All the resulting reconstructions appear to be reasonably good also. 

In the Books images zoomed in (Figure \ref{fig:Books_6400_25600_tv_tikh_zoom}, we see how the selected TV parameter is too large: we lose some of the small details. This shows that it can be tricky to choose the threshold based on the dual-grid SSIM-curve if the shape is not a nice sigmoid.

\section*{Acknowledgements}

We thank Ensio Suonperä for support in the implementation of the bilevel optimization method for parameter selection comparison.\\

\noindent This work was supported by the Research Council of Finland (Flagship of Advanced Mathematics for Sensing, Imaging and Modelling grant 359182).

\bibliographystyle{abbrv}
\bibliography{mybibliography}

@article{niinimaki2016multiresolution,
  title={Multiresolution parameter choice method for total variation regularized tomography},
  author={Niinimaki, Kati and Lassas, Matti and Hamalainen, Keijo and Kallonen, Aki and Kolehmainen, Ville and Niemi, Esa and Siltanen, Samuli},
  journal={SIAM journal on imaging sciences},
  volume={9},
  number={3},
  pages={938--974},
  year={2016},
  publisher={SIAM}
}

@book{hansen2006deblurring,
  title={Deblurring images: matrices, spectra, and filtering},
  author={Hansen, Per Christian and Nagy, James G and O'leary, Dianne P},
  year={2006},
  publisher={SIAM}
}

@book{engl1996regularization,
  title={Regularization of inverse problems},
  author={Engl, Heinz Werner and Hanke, Martin and Neubauer, Andreas},
  volume={375},
  year={1996},
  publisher={Springer Science \& Business Media}
}

@article{davoli2024dyadic,
  title={Dyadic partition-based training schemes for TV/TGV denoising},
  author={Davoli, Elisa and Ferreira, Rita and Fonseca, Irene and Iglesias, Jos{\'e} A},
  journal={Journal of Mathematical Imaging and Vision},
  pages={1--39},
  year={2024},
  publisher={Springer}
}

@ARTICLE{SSIM,
  author={Zhou Wang and Bovik, A.C. and Sheikh, H.R. and Simoncelli, E.P.},
  journal={IEEE Transactions on Image Processing}, 
  title={Image quality assessment: from error visibility to structural similarity}, 
  year={2004},
  volume={13},
  number={4},
  pages={600-612},
  doi={10.1109/TIP.2003.819861}}

@article{purisha2017controlled,
  title={Controlled wavelet domain sparsity for X-ray tomography},
  author={Purisha, Zenith and Rimpel{\"a}inen, Juho and Bubba, Tatiana and Siltanen, Samuli},
  journal={Measurement Science and Technology},
  volume={29},
  number={1},
  pages={014002},
  year={2017},
  publisher={IOP Publishing}
}

@article{chambolle2011first,
  title={A first-order primal-dual algorithm for convex problems with applications to imaging},
  author={Chambolle, Antonin and Pock, Thomas},
  journal={Journal of mathematical imaging and vision},
  volume={40},
  pages={120--145},
  year={2011},
  publisher={Springer}
}

@article{meaney2024image,
  title={Image Reconstruction in Cone Beam Computed Tomography Using Controlled Gradient Sparsity},
  author={Meaney, Alexander and Brix, Mikael AK and Nieminen, Miika T and Siltanen, Samuli},
  journal={arXiv preprint arXiv:2412.07465},
  year={2024}
}

@article{purisha2018automatic,
  title={An automatic regularization method: An application for 3-D X-ray micro-CT reconstruction using sparse data},
  author={Purisha, Zenith and Karhula, Sakari S and Ketola, Juuso H and Rimpel{\"a}inen, Juho and Nieminen, Miika T and Saarakkala, Simo and Kr{\"o}ger, Heikki and Siltanen, Samuli},
  journal={IEEE transactions on medical imaging},
  volume={38},
  number={2},
  pages={417--425},
  year={2018},
  publisher={IEEE}
}

@article{hamalainen2013sparse,
  title={Sparse tomography},
  author={Hamalainen, Keijo and Kallonen, Aki and Kolehmainen, Ville and Lassas, Matti and Niinimaki, Kati and Siltanen, Samuli},
  journal={SIAM Journal on Scientific Computing},
  volume={35},
  number={3},
  pages={B644--B665},
  year={2013},
  publisher={SIAM}
}

@article{anzengruber2009morozov,
  title={Morozov's discrepancy principle for Tikhonov-type functionals with nonlinear operators},
  author={Anzengruber, Stephan W and Ramlau, Ronny},
  journal={Inverse Problems},
  volume={26},
  number={2},
  pages={025001},
  year={2009},
  publisher={IOP Publishing}
}

@article{lukas2008strong,
  title={Strong robust generalized cross-validation for choosing the regularization parameter},
  author={Lukas, Mark A},
  journal={Inverse Problems},
  volume={24},
  number={3},
  pages={034006},
  year={2008},
  publisher={IOP Publishing}
}

@article{wen2018using,
  title={Using generalized cross validation to select regularization parameter for total variation regularization problems},
  author={Wen, You-Wei and Chan, Raymond Honfu},
  journal={Inverse Problems and Imaging},
  volume={12},
  number={5},
  pages={1103--1120},
  year={2018},
  publisher={AIMS}
}

@article{hansen1992analysis,
  title={Analysis of discrete ill-posed problems by means of the L-curve},
  author={Hansen, Per Christian},
  journal={SIAM review},
  volume={34},
  number={4},
  pages={561--580},
  year={1992},
  publisher={SIAM}
}

@article{clason2010duality,
  title={A duality-based splitting method for $\ell^{1}$-TV image restoration with automatic regularization parameter choice},
  author={Clason, Christian and Jin, Bangti and Kunisch, Karl},
  journal={SIAM Journal on Scientific Computing},
  volume={32},
  number={3},
  pages={1484--1505},
  year={2010},
  publisher={SIAM}
}

@article{dong2011automated,
  title={Automated regularization parameter selection in multi-scale total variation models for image restoration},
  author={Dong, Yiqiu and Hinterm{\"u}ller, Michael and Rincon-Camacho, M Monserrat},
  journal={Journal of Mathematical Imaging and Vision},
  volume={40},
  pages={82--104},
  year={2011},
  publisher={Springer}
}

@article{wen2011parameter,
  title={Parameter selection for total-variation-based image restoration using discrepancy principle},
  author={Wen, You-Wei and Chan, Raymond H},
  journal={IEEE Transactions on Image Processing},
  volume={21},
  number={4},
  pages={1770--1781},
  year={2011},
  publisher={IEEE}
}

@article{toma2015iterative,
  title={Iterative choice of the optimal regularization parameter in TV image restoration},
  author={Toma, Alina and Sixou, Bruno and Peyrin, Fran{\c{c}}oise},
  journal={Inverse Probl. Imaging},
  volume={9},
  number={4},
  pages={1171--1191},
  year={2015}
}

@article{chen2014automatic,
  title={An automatic regularization parameter selection algorithm in the total variation model for image deblurring},
  author={Chen, Ke and Piccolomini, E Loli and Zama, Fabiana},
  journal={Numerical Algorithms},
  volume={67},
  pages={73--92},
  year={2014},
  publisher={Springer}
}

@article{kindermann2014numerical,
  title={A numerical study of heuristic parameter choice rules for total variation regularization},
  author={Kindermann, Stefan and Mutimbu, Lawrence D and Resmerita, Elena},
  journal={Journal of Inverse and Ill-Posed Problems},
  volume={22},
  number={1},
  pages={63--94},
  year={2014},
  publisher={Walter de Gruyter GmbH}
}

@article{pragliola2023admm,
  title={ADMM-based residual whiteness principle for automatic parameter selection in single image super-resolution problems},
  author={Pragliola, Monica and Calatroni, Luca and Lanza, Alessandro and Sgallari, Fiorella},
  journal={Journal of Mathematical Imaging and Vision},
  volume={65},
  number={1},
  pages={99--123},
  year={2023},
  publisher={Springer}
}

@article{Tik63,
    author = {Tikhonov, A.N.},
    title = {On the solution of ill-posed problems and the method of regularization} ,
    journal = {Dokl. Akad. Nauk SSSR} ,
    year = {1963} 
}

@article{doi:10.1137/080724265,
author = {Wang, Yilun and Yang, Junfeng and Yin, Wotao and Zhang, Yin},
title = {A New Alternating Minimization Algorithm for Total Variation Image Reconstruction},
journal = {SIAM Journal on Imaging Sciences},
volume = {1},
number = {3},
pages = {248-272},
year = {2008},
doi = {10.1137/080724265},
}

@article{doi:10.1137/040605412,
author = {Osher, Stanley and Burger, Martin and Goldfarb, Donald and Xu, Jinjun and Yin, Wotao},
title = {An Iterative Regularization Method for Total Variation-Based Image Restoration},
journal = {Multiscale Modeling \& Simulation},
volume = {4},
number = {2},
pages = {460-489},
year = {2005},
doi = {10.1137/040605412},
}

@ARTICLE{VogelTV98,
  author={Vogel, C.R. and Oman, M.E.},
  journal={IEEE Transactions on Image Processing}, 
  title={Fast, robust total variation-based reconstruction of noisy, blurred images}, 
  year={1998},
  volume={7},
  number={6},
  pages={813-824},
  doi={10.1109/83.679423}}

@INPROCEEDINGS{RudinOsher94,
  author={Rudin, L.I. and Osher, S.},
  booktitle={Proceedings of 1st International Conference on Image Processing}, 
  title={Total variation based image restoration with free local constraints}, 
  year={1994},
  volume={1},
  number={},
  pages={31-35 vol.1},
  doi={10.1109/ICIP.1994.413269}}

@article{harrach2020beyond,
  title={Beyond the Bakushinkii veto: regularising linear inverse problems without knowing the noise distribution},
  author={Harrach, Bastian and Jahn, Tim and Potthast, Roland},
  journal={Numerische Mathematik},
  volume={145},
  number={3},
  pages={581--603},
  year={2020},
  publisher={Springer}
}

@article{bauer2008regularization,
  title={Regularization independent of the noise level: an analysis of quasi-optimality},
  author={Bauer, Frank and Rei{\ss}, Markus},
  journal={Inverse Problems},
  volume={24},
  number={5},
  pages={055009},
  year={2008},
  publisher={IOP Publishing}
}

@article{neubauer2008convergence,
  title={The convergence of a new heuristic parameter selection criterion for general regularization methods},
  author={Neubauer, Andreas},
  journal={Inverse Problems},
  volume={24},
  number={5},
  pages={055005},
  year={2008},
  publisher={IOP Publishing}
}

@article{bakushinskii1984remarks,
  title={Remarks on choosing a regularization parameter using the quasi-optimality and ratio criterion},
  author={Bakushinskii, AB},
  journal={USSR Computational Mathematics and Mathematical Physics},
  volume={24},
  number={4},
  pages={181--182},
  year={1984},
  publisher={Elsevier}
}

@article{suonpera2024linearly,
  title={Linearly convergent bilevel optimization with single-step inner methods},
  author={Suonper{\"a}, Ensio and Valkonen, Tuomo},
  journal={Computational Optimization and Applications},
  volume={87},
  number={2},
  pages={571--610},
  year={2024},
  publisher={Springer}
}

@misc{suonpera_2023_7974062,
  author       = {Ensio Suonperä},
  title        = {Codes for "Linearly convergent bilevel
                   optimization with single-step inner methods"
                  },
  month        = may,
  year         = 2023,
  publisher    = {Zenodo},
  version      = 2,
  doi          = {10.5281/zenodo.7974062},
  url          = {https://doi.org/10.5281/zenodo.7974062},
  note         =  {https://doi.org/10.5281/zenodo.7974062},
}

@article{kunisch2013bilevel,
  author  = {Kunisch, Karl and Pock, Thomas},
  title   = {A Bilevel Optimization Approach for Parameter Learning in Variational Models},
  journal = {SIAM Journal on Imaging Sciences},
  year    = {2013},
  volume  = {6},
  number  = {2},
  pages   = {938--983},
  doi     = {10.1137/120882706}
}

@article{delosreyes2017bilevel,
  author  = {De los Reyes, Juan Carlos and Sch{\"o}nlieb, Carola-Bibiane and Valkonen, Tuomo},
  title   = {Bilevel Parameter Learning for Higher-Order Total Variation Regularisation Models},
  journal = {Journal of Mathematical Imaging and Vision},
  year    = {2017},
  volume  = {57},
  number  = {1},
  pages   = {1--25},
  doi     = {10.1007/s10851-016-0662-8}
}

@article{arridge2019solving,
  author  = {Arridge, Simon and Maass, Peter and {\"O}ktem, Ozan and Sch{\"o}nlieb, Carola-Bibiane},
  title   = {Solving Inverse Problems Using Data-Driven Models},
  journal = {Acta Numerica},
  year    = {2019},
  volume  = {28},
  pages   = {1--174},
  doi     = {10.1017/S0962492919000059}
}

\end{document}